\documentstyle[12pt,subeqnarray,amssymb]{article}%

\newcommand{\al}{\alpha }
\def\b{\beta }
\newcommand{\g}{\gamma }        
\def\d{\delta }

\newcommand{\th}{\theta }       \newcommand{\Th}{\Theta }
\newcommand{\k}{\kappa }
\def\l{\lambda }            
\newcommand{\m}{\mu }
\newcommand{\n}{\nu }
\newcommand{\x}{\xi }         
\newcommand{\p}{\pi }         

\newcommand{\s}{\sigma }        

\newcommand{\f}{\phi }         
\def\ph{\varphi }
\newcommand{\ch}{\chi }
\newcommand{\ps}{\psi }        \newcommand{\Ps}{\Psi } 
\newcommand{\om}{\omega }       \newcommand{\Om}{\Omega }

\def\R{{\Bbb R} }     
\def\C{{\Bbb C} }
\def\T{{\Bbb T} }
\def\Z{{\Bbb Z} }
\def\1{{\bold 1} }  
\def\K{{\Bbb K} }

 \def\vs #1{\vspace*{#1cm} }         \def\hs #1{\hspace*{#1cm} }
      \newcommand{\ovl}{\overline}
\newcommand{\frs}{\frenchspacing\ }

 \def\lb #1{\label{#1} }
 \def\slb #1{\slabel{#1} }
 \def\rep #1{(\ref{#1})}
\newcommand{\bib}{\bibitem}
\newcommand{\bit}{\begin{itemize} }
\newcommand{\eit}{\end{itemize}}
\newcommand{\bea}{\begin{eqnarray}}   \newcommand{\bean}{\begin{eqnarray*}}
\newcommand{\eea}{\end{eqnarray}}     \newcommand{\eean}{\end{eqnarray*}}
\newcommand{\bse}{\begin{subeqnarray}} 
\newcommand{\ese}{\end{subeqnarray}}
\newcommand{\non}{\nonumber}
\def\||{\Vert }
\def\cop{\triangle }
\newcommand{\fal}{\forall }
\newcommand{\apli}{\rightarrow }    \newcommand{\map}{\mapsto }

\newcommand{\bra}{\langle }        \newcommand{\ket}{\rangle }
\def\pd{\otimes }             
\def\+{\oplus }
\newcommand{\pex}{\wedge }
 \def\cali #1{{\cal #1} }
\def\Tr{\mbox{Tr} }
 \def\del #1{\partial_{#1} }
\def\tc{\circledast }
\def\Lo{L_{\Omega} }
\def\KsH{{\Bbb K}^s(H_3) }
\def\ko{\kappa^{\Omega} }
\def\kt{\kappa^{\Theta} }
\def\copo{\triangle^{\Omega} }
\def\copt{\triangle^{\Theta} }
\def\pho{\varphi^{\Omega} }
\def\pht{\varphi^{\Theta} }
\def\omo{\omega^{\Omega} }
\def\omt{\omega^{\Theta} }
\def\oml{\hat\omega^{\n} }

\begin{document}

\title{\bf Projective Fourier Duality\\ and Weyl
Quantization\thanks{Work supported by CAPES and CNPq, Brazil.}}

\author{R. Aldrovandi\thanks{e-mail: ra@axp.ift.unesp.br}\ \ and L.A.
Saeger\thanks{e-mail: saeger@qcd.th.u-psud.fr}\\ Instituto de F\'\i sica
Te\'orica - UNESP\\ Rua Pamplona 145\\ 01405-900 - S\~ao Paulo, SP\\ Brazil}
\date{}
\maketitle

The Weyl-Wigner correspondence prescription, which makes large use of Fourier
duality, is reexamined from the point of view of Kac algebras, the most general background for noncommutative Fourier analysis allowing for that property. It is shown how the standard Kac structure has to be extended in order to accommodate the physical requirements. An Abelian and a symmetric {\em projective Kac algebras} are shown to provide, in close parallel to the
standard case, a new dual framework and a well-defined notion of {\em
projective Fourier duality} for the group of translations on the plane. The
Weyl formula arises naturally as an irreducible component of the duality
mapping between these projective algebras.

\vs{.7}

{\bf Keywords}: Weyl Quantization, Weyl-Wigner Correspondence, Projective
Group Duality, Projective Kac Algebras, Noncommutative Harmonic Analysis.

\vs{.5}

Internat. J. Theor. Phys. {\bf 36}(3), 573-612 (1997)\hspace*{\fill} funct-an/9608004

\newpage

\section{Introduction}

In its broadest meaning, the word ``quantization'' signifies the passage from
the classical to the quantum description of a system. The most complete
classical description being found in the Hamiltonian formalism,
the natural path to take is quantization on phase space. This is the main
appeal of the Weyl-Wigner approach, which realizes the correspondence
principle by attributing a quantum operator to each classical dynamical
variable via a Fourier transformation of its density. Conversely, it also
attributes a c-number function to each operator by another Fourier
transformation, this time involving an integration on operator space. All
this supposes the possibility of performing two-way Fourier transformations,
that is, of doing a transformation {\em and its inverse}. Two points should
be retained: (i) the integration over operator space, as usually presented,
is purely formal and should be better defined; (ii) the two-way Fourier
transformations make use of a deep property of harmonic analysis, Fourier
duality. This property only holds under very severe conditions. Actually, at
least as usually presented, the Weyl-Wigner formalism supposes a very
particular kind of that duality, the Pontryagin duality, which should not be
expected to be at work when the phase space is not, for each degree of
freedom, the plane $\R^2$ which models vector phase spaces. Pontryagin
duality is valid only when the group of linear symplectomorphisms
(transformations preserving the phase space symplectic structure) is Abelian.
This group is, for $\R^2$, the translation group $\R^2$. To quantize on more
general phase spaces, we are bound to consider the general approach to
Fourier transformations, which requires Kac algebras. We have shown elsewhere
\cite{rasa} how this general formalism, in all its complexity, is necessary
even for the simplest non-trivial phase space, the half-plane. We intend here
to revisit the apparently well-known $\R^2$ case from this point of view.
Because of its vector space structure, the plane would seem to dispense with
a more involved treatment. We shall see that this is not so. It actually
conceals a lot of structure under the appearance of simplicity and, due to
its non-trivial cohomology, requires an extension of the very concept of Kac
algebra. Furthermore, through the pioneering work of Ref.~\cite{segal} and the subsequent introduction of weights in the sixties, the general
approach provides a precise meaning to the otherwise mysterious integration
over operator space.
 
Kac algebras are the most general structures known nowadays on which Fourier
analysis can be realized in its integrity. Their rather involved axioms are
essential to the most demanding of the properties attached to harmonic
analysis, precisely the duality above mentioned. As soon as we depart from
the case of functions on Abelian groups, for which the Pontryagin
group-to-group duality holds in all its simplicity, Fourier transforms and
their inverses can be defined only for functions on domain spaces much more
sophisticated than groups. Starting from functions on groups which are
separable and locally compact, we arrive necessarily to Kac algebras, which
are Hopf-von Neumann algebras endowed with Haar weights. This means that they
are noncommutative spaces on which we known how to perform (noncommutative)
integration. Roughly speaking, whenever we do Fourier analysis, we are
supposing the presence (explicit or not) of Kac algebras. We propose here to
bring to light the algebras behind the apparently simple case of the plane
$\R^2$. And here comes an important point. As they are known today, Kac
algebras are related to linear representations and as such they are not
sophisticated enough to cope with the problem. In order to apply to Quantum
Mechanics, the Kac structure may require an extension to projective
representations, and this is precisely what happens in the usual Weyl-Wigner
formalism. The situation is rather curious. On one hand, so much is
``degenerated'' in this simplest of all cases (the space dual to $\R^2$ is
$\R^2$ itself, which coincides also with the group manifold of linear
symplectomorphisms) that we get the impression that the intricacies of the
general formalism can be overlooked; on the other, because of its nontrivial
cohomology, it requires an extension to projective representations, which is
not necessary in other, more complicated, situations. For example, no
extension is required when the phase space is the half-plane \cite{rasa}. The
central extension of $\R^2$ is, roughly speaking, the Heisenberg group $H_3$.
No extension is a particular case: it can be seen as a trivial extension. We
can say that extensions, trivial or not, are required in the generic case and
their effects on the standard structures have to be studied.

The essential notation is introduced in sections~\ref{cqmph} and
\ref{heigr}, which sum up the usual lore on quantization on phase space and the
Heisenberg group. Projective representations of a group can be obtained from
the linear representations of its extension. We arrive thus at the projective
representations of $\R^2$ from the linear representations of $H_3$. For Kac
algebras, a parallel procedure will be used: we start from the
well-established Kac algebra duality for $H_3$ and then proceed to find the
projective Kac algebras of $\R^2$. Actually, a pair of Kac algebras is
necessary to the materialization of duality. One, called the Abelian Kac
algebra, has for elements the $L^{\infty}$-functions. The other, the
symmetric Kac algebra, includes the left-regular representations. The problem
lies in the fact that the symmetric Kac algebra for the Heisenberg group
$H_3$ is generated by (linear !) left-regular representations while the Weyl
kernels are irreducible projective operators. This is reviewed in
section~\ref{fdhg}. To go from the $H_3$-Kac duality to the desired
projective algebras, two steps are involved: projection and decomposition
into irreducibles, in this order or in its inverse. The extensions so
obtained are far from trivial. Kac algebras are, to begin with, Hopf-von
Neumann algebras, and the necessary extensions involve generalizations of
some of the current concepts on Hopf algebras. Though most of the axioms
remain unchanged, some of the usual requirements, valid for linear
representations, must be extended to their projective counterparts. This
generalization to projective Kac algebras is presented in section~\ref{prka}.
It leads to the general notion of projective Kac algebra. An extended duality
comes out, a projective Kac duality leading to a projective Fourier duality.
It is necessary to introduce some notions like projective coinvolution,
coprojective coinvolution, as well as to extend the usual axioms concerning
anti-(co)automorphisms. In section~\ref{wqd}, the Weyl-Wigner correspondence
is recast into the Fourier duality language. Weyl's formula comes up as an
irreducible component of the duality mapping between the previously obtained
projective algebras. A tentative prescription for quantization on general
phase spaces is sketched in the final considerations.

\section{Classical and Quantum Mechanics in Phase Space}\lb{cqmph}

The classical picture of Mechanics, as is well known, can be described
geometrically in terms of the symplectic structure of the phase space
\cite{arn}. The simplest symplectic manifold is just $\R^2$, the most usual
arena of classical dynamical systems with one degree of freedom. On this
phase space, with its trivial symplectic structure $\om=dp\pex dq$, acts by
symplectomorphisms the Abelian group -- also denoted $\R^2$ -- of
two-dimensional translations. To classical Hamiltonian dynamical systems on
$\R^2$, associated to Hamiltonian functions $H$, correspond symplectic 
Hamiltonian vector fields $X_H$ satisfying the Hamilton equations
\bea
 i_{X_H}\om=-dH.
\lb{hameq} 
\eea
The non-degeneracy of the symplectic form implies a local isomorphism between
vector fields and 1-forms (see \rep{hameq}), and a homomorphism between
vector fields and $C^{\infty}$-functions. Vector fields constitute a Lie
algebra by the Lie bracket, whose isomorphic image on the space of functions
is the Poisson bracket
\bean
 \{f,g\}=-\om(X_f,X_g).
\eean 
According to Dirac \cite{dir}, in order to quantize on such a space, we must
be sure that there exists a faithful correspondence between this Poisson
algebra and an operatorial algebra. The closest operator algebra we have at
hand is the Lie algebra of the group acting on the phase space by
symplectomorphisms. In the Euclidean case, since the translation group is
Abelian, we must central-extend it to the Heisenberg group in order to have
the isomorphism of the group and the Poisson Lie algebras. Such an
isomorphism allows the construction of a faithful quantization map on this
phase space \cite{ish}.

From the point of view of Harmonic Analysis, the Pontryagin duality for the
Abelian group $\R^2$ ensures that the Fourier transform and its inverse
constitute an isomorphism between the Abelian convolution algebra $L^1(\R^2)$
and the also Abelian algebra $L^{\infty}(\R^2)$ of essentially bounded
functions with pointwise product, both contained in $C^{\infty}(\R^2)$. The
Fourier transform of an $L^1$-function $f$ is the $L^{\infty}$-function
\bea
[\cali Ff](x,y)=\tilde f(x,y)=\frac 1{2\p}\int_{\R^2}dqdp\,
f(p,q)\,e^{-i(yq+xp)},
\lb{fotr} 
\eea
where the kernel $e^{i(yq+xp)}\equiv\ch_{(x,y)}(q,p)$ is a character
(one-dimensional irreducible representation) of the group $\R^2$. Characters
satisfy the orthogonality-completeness relation
\bean
\int_{\R^2}
dxdy\,\ch_{(x,y)}(q,p)\ovl{\ch_{(x,y)}(q',p')}=(2\p)^2\d(q-q')\d(p-p'),
\eean
which can be used to invert \rep{fotr} and write the inverse Fourier
transform as
\bea
 f(p,q)=\frac 1{2\p}\int_{\R^2}dxdy\, \tilde f(x,y)\,e^{i(yq+xp)}. 
\lb{iftr}
\eea
From this point of view, the correct way to regard these formulas is to see
the first as the Fourier transform from $L^1(\R^2)$ into
$L^{\infty}(\hat{\R^2})$, which is the same as $L^{\infty}(\R^2)$ since this
group is self-dual, $\hat{\R}^2=\R^2$, and the second as the transform from
$L^1(\hat{\R^2})$ into $L^{\infty}(\R^2)$. This is so because the Fourier
transform is defined as a mapping between the $L^1$-space of an Abelian group
into the $L^{\infty}$-space of its dual, the space of characters which is
also a group. Because $\R^2$ is self-dual, the Fourier transform turns out to
be an algebra isomorphism, mapping convolution to pointwise product
\cite{rei}. This enables us to regard formula \rep{iftr} as the inversion of
\rep{fotr}, or as the Fourier transform for $\hat R^2$. Thus, from the point
of view of Harmonic Analysis, this group is highly ``degenerate''. Extending
the domain of the transform $\cali F$ from $L^1$ to the space $\cali
S'(\R^2)$ of tempered distributions on $\R^2$, $\cali F:\cali
S'(\R^2)\apli\cali S'(\R^2)$ turns out to be a topological isomorphism
\cite{sugi,3sen}. The same happens when $\cali F$ is restricted to the space
$\cali S(\R^2)$ of rapidly decreasing functions, confirming the degeneracy
alluded to.

To go from this classical approach to a quantum picture, Weyl proposed
\cite{wey} to modify the Fourier transform formula by changing its
scalar kernel into an operatorial kernel. He wrote
\bea
\hat f_{\hbar} =\int_{\R^2}dqdp\, f(q,p)\,e^{-\frac i{\hbar}(p\hat q+q\hat p)},
\lb{wefor}
\eea
instead of $\tilde f$, where $\hat q,\ \hat p$ are the usual coordinate and
momentum operators of Euclidean Quantum Mechanics. These operators satisfy the
Heisenberg commutation relation $[\hat q,\hat p]=i\hbar$. By the Glauber
identity, the operatorial kernel $S'_{\hbar}(x,y)\equiv e^{-\frac
i{\hbar}(y\hat q+x\hat p)}$ can be written also as
\bean
 S_{\hbar}(x,y)=e^{\frac i{2\hbar}xy}\,U(y)V(x),
\eean
in terms of the Weyl operators $U(y)=e^{-\frac i{\hbar}y\hat q}$ and
$V(x)=e^{-\frac i{\hbar}x\hat p}$. These satisfy the Weyl commutation
relation (we use the representation theory convention, by which the leftmost
operator acts first)
\bean
 U(y)V(x)=e^{-\frac i{\hbar}xy}\,V(x)U(y),
\eean
while for $S_{\hbar}$ holds
\bea
S_{\hbar}(x,y)S_{\hbar}(x',y')=e^{\frac
i{2\hbar}(xy'-yx')}S_{\hbar}(x+x',y+y').
\lb{popr}
\eea
These are consequently not linear, but projective operators. They realize an
operatorial representation of $\R^2$, otherwise impossible for such an
Abelian group. Summing up, as a realization of the classical-quantum
correspondence principle, Weyl proposed with \rep{wefor} to consider a
passage from the usual scalar Fourier transform to an operatorial one written
in terms of projective representations. The inverse way \cite{wig}, leading
from the quantum to the classical picture, involves an {\em integration} on
operator space (taking of the trace):
\bea
 f(q,p)=\Tr [S_{\hbar}^{\dagger}(q,p)\hat f_{\hbar}].
\lb{inwef}
\eea
Since the projective operator product \rep{popr} carries a natural twisting,
and since formulas \rep{wefor}, \rep{inwef} ought to represent an algebra
isomorphism, the corresponding convolution in $L^1(\R^2)$ also gets twisted:
\bean
 \hat f_{\hbar}\cdot\hat g_{\hbar}&=&\int_{\R^2\times\R^2}dxdydx'dy'
 f(x,y)g(x',y')\,e^{\frac i{2\hbar}(xy'-yx')}S_{\hbar}(x+x',y+y')\non\\
&=&\int_{\R^2}dx''dy''\,(f\tc g)(x'',y'')\,S_{\hbar}(x'',y''),
\eean
where
\bean
(f\tc g)(x'',y'')=\int_{\R^2} dxdy\,e^{\frac i{2\hbar}(xy''-yx'')}
f(x,y)g(x''-x,y''-y).
\eean
Formulas \rep{wefor} and \rep{inwef} provide a two-way correspondence between
the classical ($L^1$-functions) and the quantum pictures. The further
correspondence to $L^{\infty}$-functions is provided by the Fourier transform.
The Fourier transform of a twisted convolution of two functions gives rise to
the twisted (noncommutative) product of their Fourier transforms, which
characterizes a deformation of the Abelian algebra of the pointwise product.
This two-way classical-quantum procedure is the {\em Weyl-Wigner correspondence
prescription}.

Projective representations of a group are generally obtained from the linear
representations of its central extension \cite{barg}, in our case the
Heisenberg group. Harmonic Analysis on general locally compact groups like
the Heisenberg group is not as trivial as that on the Abelian ones. Because
such groups have infinite dimensional irreducible representations, finite
algebras have not enough structure to host a duality. We must deal with
semifinite von Neumann algebras endowed with additional structure, the Kac
algebras. From the relation between the irreducible representations of the
Heisenberg group and the projective operators appearing in the Weyl-Wigner
formalism, we can relate Harmonic Analysis on the Heisenberg group to a
projective Harmonic Analysis on $\R^2$ and ``explain'' the origin of the
Weyl-Wigner formulas. This will be done in the last sections, after we have
established some facts on the Heisenberg group in section
\ref{heigr}, and reviewed the duality theory for it in terms of Kac algebras
in section \ref{fdhg}.

\section{The Heisenberg Group}\lb{heigr}

In this work the three-dimensional Heisenberg group $H_3$ is regarded as the
central extension of the two-dimensional Abelian group of translations on the
plane by the torus $\T$. We shall use the notation
$(x,\al)=(x_1,x_2,e^{i\th}),\ x_1,x_2,\in\R$, $\th\in\R/{2\pi}$, to denote the
elements and coordinates of $H_3$. As is well known, the second cohomology
space $H^2(\R^2,\R/2\p)$ of cocycles from $\R^2$ to $\R$ (mod $2\p$) is not
trivial \cite{tuy}. Since 2-cocycles classify central extensions, inequivalent
2-cocycles give rise to inequivalent central extensions. Thus, for a chosen 
cocycle $\Om\in H^2(\R^2,\R/2\p)$, e.g. 
\bea
 \Om(x,y)=\frac 12(x_1y_2-y_1x_2),
\lb{tcoh}
\eea
the product on $H_3=\R^2\times\T$ is given by
\bean
 (x,\al)(y,\b)=(x+y,\al\b\,e^{i\Om(x,y)}),
\eean
where associativity is ensured by the closeness of $\Om$ in $H^2$, namely,
\bean
 \d\Om(x,y,z)=\Om(y,z)-\Om(x+y,z)+\Om(x,y+z)-\Om(x,y)=0.
\eean
The identity in $H_3$ is $(0,1)$ and the inverse element of $(x,\al)$ is
$(-x,\al^{-1})$. The following useful properties of $\Om$ are obvious from
\rep{tcoh}: $\Om(x,0)=0$, $\Om(x,y)=-\Om(y,x)$, $\Om(-x,y)=-\Om(x,y)$.

The irreducible linear representations of $H_3$ can be obtained by Mackey's
induced representation method \cite{mack,tay}. Their division into inequivalent
classes is given by the Stone-von Neumann theorem, which also provides the
unitary dual space of this group. These representations are divided into
infinite-dimensional and one-dimensional ones in the dual $\hat{H_3}=(\Z
-\{0\})\cup\R^2$, according to
\bse
\lb{stonevn} &&T_{\n}(x,\al)=e^{i\n\th}e^{\frac i2\n x_1x_2}e^{-i\n x_2\hat
q}e^{-ix_1\hat p},\hs 1\n\in\Z-\{0\}\slb{idirh}\\
&&T_{ab}(x,\al)=e^{iax_2}e^{ibx_1},\hs 1 (a,b)\in\R^2,\slb{irchar}
\ese
where the self-adjoint operators $\hat q,\ \hat p$ act on $L^2(\R)$ by
\bean
 && \hat q\ps(q)=q\ps(q)\\
 && \hat p\ps(q)=-i\del q\ps(q).
\eean
We recall that the commutation relation
\bean
 [\hat q,\hat p]=i
\eean
is a realization of the Lie algebra of $H_3$ on that Hilbert space, which is
also isomorphic to the Poisson algebra generated by the coordinates $q,\ p$
plus the constant function $1$ on the Euclidean symplectic manifold $\R^2$.

\section{Fourier Duality for the Heisenberg Group}\lb{fdhg}

Our objective is to describe the Weyl-Wigner correspondence in terms of
projective Fourier duality, that is, to find out a connection between Kac
duality and the algebra generated by irreducible projective operators. Since
the projective representations of $\R^2$ are obtained from the linear
representations of $H_3$, we should start from the well-established Kac
algebra duality for this group. This is reviewed in this section. As already
said, duality requires a pair of Kac algebras: the Abelian Kac algebra,
formed with the $L^{\infty}$-functions, and the symmetric Kac algebra,
including the left-regular representations.

\subsection{The Symmetric Kac Algebra of $H_3$}\lb{skah3}

Let us begin by introducing the symmetric Kac algebra of $H_3$, $\KsH$, which
is built on the von Neumann algebra $\cali M(H_3)$ generated by the
left-regular representation operators of the group. For details on Kac
algebras, see the book~\cite{ensc} and the quick review in the first sections
of Ref.~\cite{vai}, or still Ref.~\cite{rasa}. First recall that the
left-regular representation $L$ acts on the Hilbert space $L^2(H_3)$ of
square-integrable functions on the Heisenberg group by
\bea
 [L(x,\al)f](y,\b)=f((x,\al)^{-1}(y,\b)).
\lb{alroh}
\eea
The scalar product in this space is given by
\bean
 (f|g)_{L^2(H_3)}=\int_{H_3}dxd\al\ f(x,\al)\ovl{g(x,\al)}
\eean
and the norm by $\||f\||_2^2=(f|f)_{L^2(H_3)}$, where $dxd\al$ is the
left-and-right invariant measure on $H_3$, which is unimodular. In 
this section the spaces $L^p(H_3),\ p=1,2,\infty$ will be denoted simply by
$L^p$.

$\cali M(H_3)$ is a subalgebra of $\cali B(L^2)$, the Banach algebra of all
bounded operators on the Hilbert space $L^2$. This means that the product on
$\cali M(H_3)$ is associative, there exists a unit $I=L(0,1)$ (the identity
operator), and also an involution (taking of the dagger) such that
$I^{\dagger}=I$. The norm is defined by $\||T\||=sup\{\||T\ps\||_2,\
\||\ps\||_2=1\}$, under which $\||T^{\dagger}\||=\||T\||$ ($\cali M(H_3)$ is
an involutive Banach algebra) and $\||T^{\dagger}T\||=\||T\||^2$ (it is a
$C^*$-algebra). There is also a family of seminorms defined by
$\||T\||_{w,\ps,\f}=|(\ps|T\f)_{L^2}|,\ \ps,\f\in L^2$, whose open balls
define the {\em weak topology}. $\cali M(H_3)$ is closed in this topology.

The elements of $\KsH$ are written in terms of the generators as
\bea
 \hat f=\int_{H_3}dxd\al\ f(x,\al)L(x,\al),
\lb{lrrl}
\eea
where the coefficients $f$ are functions of compact support, whose algebra
$C(H_3)$ is dense in the convolution Banach algebra $L^1$. The product can be
written in terms of the convolution of the coefficients:
\bean
\hat f\cdot\hat g&=&\int_{H_3\times
H_3}dxd\al\,dyd\b\,f(x,\al)g(y,\b)\,L((x,\al)(y,\b))\non\\
&=&\int_{H_3}dzd\g\, (f\ast g)(z,\g)\, L(z,\g),
\eean
where the convolution on $C(H_3)$ is written
\bean
 (f\ast g)(z,\g)=\int_{H_3}dxd\al\,f(x,\al)g((x,\al)^{-1}(z,\g)).
\eean 
Formula \rep{lrrl} can also be regarded as expressing the left-regular
representation of $L^1$ induced by the left-regular representation of the
group, and as such is denoted $L(f)$.

With such elements, $\KsH$ has a structure given by the following operations:
\bit
\item a {\em product} given by the group multiplication, 
\bea
 L(x,\al)L(y,\b)=L(x+y,\al\b\,e^{i\Om(x,y)});
\lb{kspr}
\eea
\item a symmetric (wherefrom the name of this algebra) {\em coproduct},
\bea
 \hat{\cop} L(x,\al)=L(x,\al)\pd L(x,\al);
\lb{copd}
\eea
\item a {\em coinvolution},
\bea
 \hat{\k}(L(x,\al))=L^{\dagger}(x,\al)=L((x,\al)^{-1});
\lb{kscoi}
\eea
\item and a normal, faithful and semifinite (n.f.s.) Haar {\em trace},
\bea
 \hat\ph(T)=\left\{
\begin{array}{cl}
 \||f\||^2_2 & \mbox{if}\ T=\hat f^{\dagger}\cdot\hat f\\
+\infty & \mbox{otherwise}
\end{array}
\right.\hs 1 T\in\cali M(H_3)^+.
\lb{twwt}
\eea
\eit
Normal, faithful and semifinite mean respectively that: $\hat\ph(T)$ is the
upper bound of the sequence $\{\hat\ph(T_i)\}$ if $T\in\cali M(H_3)^+$ is the
upper bound of the sequence $\{T_i\}$; $\hat\ph(T)=0$ implies $T=0$, $T\in
\cali M(H_3)^+$; the algebra $span\{T\in \cali M(H_3)^+\ |\
\hat\ph(T)<\infty\}$ is $\s$-weakly dense in $\cali M(H_3)$. The $\s$-weak
topology is defined by the open balls of the family of seminorms
$\||F\||_{\s,\f_i,\ps_i}=\sum_i|(\f_i|F\ps_i)|$, where
$\sum_i\||\f_i\||^2<\infty$, $\sum_i\||\ps_i\||^2<\infty$. $\cali M(H_3)^+$
is the set of positive elements of $\cali M(H_3)$, that is, the set of
operators with positive spectrum.

Equation \rep{twwt} is coherent with 
\bean
 \hat\ph(\hat f)=f(0,1), 
\eean
for in this case $\hat\ph(\hat f^{\dagger}\cdot\hat f)=(f^*\ast f)(0,1)=
\||f\||^2_2$, where $^*$ denotes the involution on $L^1$ given by
\bea
 f^*(x,\al)=\ovl{f((x,\al)^{-1})}.
\lb{l1inv}
\eea 

The coproduct has a canonical implementation on $\cali M(H_3)$ in terms of
a unitary operator $\hat W\in\cali B(L^2)\pd L^{\infty}$, 
\bea
 \hat\cop L(x,\al)=\hat W(I\pd L(x,\al))\hat W^*.
\lb{cioc}
\eea
This {\em fundamental operator} is unique and is fixed by
\bea
 [\hat W F](x,\al;y,\b)=F((y,\b)^{-1}(x,\al);(y,\b)),
\lb{dfu}
\eea
where $F\in C(H_3\times H_3)$. Its adjoint $\hat W^*$ is given by $[\hat W^*
F](x,\al;y,\b)=F((y,\b)(x,\al);(y,\b))$. The importance of $\hat W$ and its
dual $W=\s\circ\hat W^*\circ\s$ lies in that they {\em generate the
Kac duality}, in the sense that they are the generators of the representations
linking $\K^s(H_3)$ and its dual. As a consequence of that and of \rep{cioc},
they satisfy the {\em pentagonal relation}
\bean
(I\pd\hat W)(\s\pd I)(I\pd\hat W)(\s\pd I)(\hat W\pd I)=(\hat W\pd
I)(I\pd\hat W).
\eean
In the same way, the coinvolution has a canonical implementation in terms of
the antilinear isometry $J:L^2\apli L^2$ by
\bean
 \hat\k(L(x,\al))=J\,L^{\dagger}(x,\al)\,J.
\eean
The latter is given on $C(H_3)$ by $[Jf](x,\al)=\ovl{f(x,\al)}$, and in the
case of an unimodular group like $H_3$, it also implements the involution in
the (pre)dual algebra.

\subsection{The Abelian Kac Algebra of $H_3$}

In order to be a Kac algebra, $\KsH=(\cali M(H_3),\hat\cop,\hat\k,\hat\ph)$
must satisfy a certain set of axioms. These will be presented later in the
subsection~\ref{prksh3} on the projection process. We only anticipate that
$\KsH$, as introduced above, does satisfy them. By the time being we are
interested in duality for $H_3$. The dual of $\KsH$ is obtained as the image
of the {\em Fourier representation} $\hat\l$ of the predual of $\cali
M(H_3)$. The predual $\cali M(H_3)_*$, which is isomorphic to the {\em
Fourier algebra} $A(H_3)$ of $H_3$, is the space of all $\s$-weakly
continuous linear functionals on $\cali M(H_3)$. The representative elements
of $\cali M(H_3)$ are linear forms $\hat\om_{fg}$ on $L^2$, in terms of which
the corresponding functions in $A(H_3)$ are defined by
\bea
\hat\om_{fg}(x,\al)\equiv\bra
L^{\dagger}(x,\al),\hat\om_{fg}\ket=(f\ast\check g)(x,\al)\in A(H_3),\hs 1
f,g\in L^2,
\lb{dplf}
\eea
where $\bra L(x,\al),\hat\om_{fg}\ket\equiv(L(x,\al)f|g)_{L^2}$ by definition
of $\hat\om_{fg}$, and $\check g(x,\al)=\ovl{g((x,\al)^{-1})}$. Notice that,
by applying the Cauchy-Schwartz inequality to \rep{dplf}, we find that this
function has an upper bound, that is,
$|\hat\om_{fg}(x,\al)|\le\||f\||_2\||g\||_2<\infty$, and consequently
$A(H_3)\subset L^{\infty}$. The product in the predual $A(H_3)$ is obtained
by duality from the coproduct in $\K^s(H_3)$,
\bea
\bra L^{\dagger}(x,\al), \hat\om_{fg}\cdot\hat\om_{hl} \ket=\bra \hat\cop
L^{\dagger}(x,\al), \hat\om_{fg}\pd\hat\om_{hl} \ket,
\lb{ddp}
\eea
and, as follows trivially from \rep{copd}, is the Abelian pointwise product.
The involution $^o$ in $A(H_3)$ also follows by duality from
\bea
\bra L^{\dagger}(x,\al), \hat\om_{fg}^o \ket=\ovl{\bra
\k(L^{\dagger}(x,\al))^{\dagger}, \hat\om_{fg} \ket},
\lb{ddi}
\eea
and is simply the complex conjugation implemented by $J$.

To find out the Fourier representation $\hat\l$, defined by
\bea
 [\hat\l(\hat\om)f](x,\al)=[(\hat\om\circ\hat\k\pd id)(\hat\cop\hat
f)]_{\hat\ph}(x,\al),\hs 1 f\in L^2,
\lb{frdef}
\eea
we may use the formula
\bea
(\hat W(f\pd g)|h\pd l)_{L^2\pd L^2}=(g|\hat\l(\hat\om_{hf})l)_{L^2},\hs 1
f,g,h,l\in L^2,
\lb{fcgr}
\eea
which relates $\hat\l$ to the dual of its generator $W$. Computing the double
scalar product in \rep{fcgr}, taking into account \rep{dfu}, and identifying
$\hat\om_{hf}$ from \rep{dplf} we get $\hat\l=id$. This means that the Kac
algebra dual to $\K^s(H_3)$ is built on the von Neumann algebra $L^{\infty}$ of
measurable and essentially bounded functions on the Heisenberg group. The
coproduct, coinvolution and trace thus obtained, together with the pointwise 
product and the involution, satisfy the Kac algebra axioms. This Abelian Kac
algebra $\K^a(H_3)$ on $L^{\infty}$ is then defined by the following structure:
\bse
 [f\cdot g](x,\al)&=&f(x,\al)g(x,\al);\\
 \1&=&1,\ \mbox{such that}\ 1(x,\al)=1\ \fal (x,\al);\\
 \cop(f)((x,\al)\pd(y,\b))&=&f((x,\al)(y,\b));\slb{copi}\\
 \k(f)(x,\al)&=&f((x,\al)^{-1});\\
 \ph(f)&=&\int_{H_3}dxd\al\,f(x,\al),\hs 1 f\in L^{\infty +}.\slb{kaph}
\ese
The positive elements are the positive definite functions in
$L^{\infty}$. This von Neumann algebra is also a subalgebra of $\cali B(L^2)$,
which acts on $L^2$ by pointwise multiplication. Its norm is given
by $\||f\||_{\infty}=ess.sup.|f(x)|$, which is the smallest number $C$ $(0\le C
<\infty)$ such that $|f(x)|\le C$ locally almost everywhere \cite{rei}. The
predual of $L^{\infty}$ is just $L^1$, the convolution algebra with involution
given by \rep{l1inv}. Needless to say that its structure is also obtainable by
duality relations similar to \rep{ddp} and \rep{ddi}, but now between
$L^{\infty}$ and $L^1$. The fundamental operator for this algebra is $W$, which
implements $\cop$ and is given by
\bean
 [WF](x,\al;y,\b)=F((x,\al);(x,\al)(y,\b)),
\eean
while the dual $\hat J$ of $J$ is given by $[\hat
Jf](x,\al)=\ovl{f(-x,\al^{-1})}$, in terms of which we have $\k(f)=\hat J\,\ovl
f\,\hat J$ on $L^2$.

The duality $\K^s(H_3)-\K^a(H_3)$ for the Heisenberg group will be complete
when $L^{\infty}_*=L^1$ is represented in $\cali M(H_3)$. This is carried out
by the Fourier representation $\l$, dual of $\hat\l$, which is just the
regular representation of $L^1$ restricted to act on $L^{\infty}\cap L^2$.
Another way to see that, and in fact to deduce it, is to use the dual of
formula \rep{fcgr},
\bea
(W(f\pd g)|h\pd l)_{L^2\pd L^2}=(g|\l(\om_{hf})l)_{L^2},\hs 1 f,g,h,l\in L^2,
\lb{fcdr}
\eea
where $\om_{fg}\in L^1$ is defined by $\bra
h,\om_{fg}\ket=(hf|g)\therefore\om_{fg}=f\ovl g$. We obtain
\bean
 \l(f)=\int_{H_3}dxd\al\, f(x,\al)\,L(x,\al),\hs 1 f\in L^1,
\eean
whose image is just $\cali M(H_3)$. Recall that formula \rep{fcdr} is a
consequence of the dual $\hat W=\s\circ W^*\circ\s$ being the
generator of $\l$. To see that and to understand what the word {\em generator}
really means, define $\f\in L^2(H_3,L^2)$ by $[\f(y,\b)](x,\al)=F(x,\al;y,\b)$,
where $F\in L^2\pd L^2$. These spaces are isomorphic. Recall also that, as a
representation, the operator $L$ is a bounded map between $H_3$ and $\cali
B(L^2)$. Then for $(y,\b)$ fixed, $L(y,\b)\in\cali B(L^2)$, $\f(y,\b)\in L^2$,
and we have
\bean
[L(y,\b)\f(y,\b)](x,\al)&=&[\f(y,\b)]((y,\b)^{-1}(x,\al))\non\\
&=&F((y,\b)^{-1}(x,\al);(y,\b)),
\eean
which is just $[\hat WF](x,\al;y,\b)$ as given in \rep{dfu}. That is, $L:
H_3\apli\cali B(L^2)$, which induces (generates) $\l$, can be seen as the
operator $\hat W\in\cali B(L^2)\pd L^{\infty}$. This is put in compact form 
as $\l(f)=(id\pd f)(\hat W)$. 

As a final remark regarding such Kac algebras, notice that both $\K^a(H_3)$
and $\K^s(H_3)$ are represented on $L^2(H_3)$ by the Gelfand-Naimark-Segal
(GNS) construction. The first is represented by the inclusion of $L^{\infty}$
and the latter by the inclusion of the $L^1$-coefficients.

\section{Projective Kac Algebras}\lb{prka}

In this section we project the symmetric and Abelian Kac algebras of $H_3$
into algebras related to the projective representations of $\R^2$ and obtain a 
projective duality extension. It is worthwhile to spend some time in the 
definition of the projective representations, since they are crucial for the 
projection process.

We start from Bargmann's \cite{barg} method to obtain projective
representations of a group from the linear representations of its central
extensions. The representations of the central extension giving rise to
projective representations are those reducing to the identity when restricted
to the central subgroup. In our case we have the left-regular representations
of $H_3$, which act on $L^2(H_3)$ by \rep{alroh}. As before, we will
concentrate on the central extension defined by the cocycle $\Om$ introduced
in \rep{tcoh}. It is clear from \rep{alroh} that the restriction of $L$ to
$\T$ is not the identity representation, so that we must make $L$ to act on
another space, suitable to our purposes. By Mackey's induced representation
method \cite{mack}, the left-regular representation can be regarded as induced
by the identity representation of the subgroup $\{e\}$. If in the induction
process we change any other subgroup for $\{e\}$, the resulting representation
is called {\em quasi-regular} \cite{bar}. Thus, since $\T$ is central, its
behaviour is equivalent to that of $\{e\}$, which enables us to interpret
the representations induced by the identity representation of $\T$ as the
regular representations acting on another space. By that method,
$L(x,\al)$ should act on a Hilbert space isomorphic to $L^2(\R^2)$, which we
call here $H(H_3)$. Its elements are square-integrable
functions when restricted to $\R^2$ and, furthermore, satisfy
$f((x,\al)(0,\b))=\b^{-1}\,f(x,\al)$. Since
$(x,\al)=(x,1)(0,\al)=(0,\al)(x,1)$, we have the decomposition
\bea
 f(x,\al)=\al^{-1}\, f(x,1)\equiv \al^{-1}\, f(x),
\lb{hsprf}
\eea
where the same notation $f$ for functions on $H_3$ and on $\R^2$ is used. By
this natural {\em projection} of $L^2(H_3)$ into $L^2(\R^2)$, \rep{alroh} can
be rewritten as
\bean
 [L(x,\al)f](y)=\al\,e^{-i\Om(-x,y)}\,f(y-x),
\eean
which does reduce to the identity when restricted to $\T$, namely
$[L(0,\al)f](y)=\al\,f(y)$. The respective projective representation of
$\R^2$ is then defined on $L^2(\R^2)$ by
\bea
 [\Lo(x)f](y)\equiv [L(x,1)f](y)=e^{i\Om(x,y)}\,f(y-x).
\lb{aclo}
\eea
From what has been said above we can also write the decomposition
of $L(x,\al)$ as
\bea
 L(x,\al)=\al\,\Lo(x).
\lb{dlrip}
\eea
As a consequence, the $\R^2$ operation (sum) is now represented by
\bea
 \Lo(x)\Lo(y)=e^{i\Om(x,y)}\Lo(x+y),
\lb{prirp}
\eea
which characterizes a projective representation.

\subsection{Projective Kac Algebras of the Translation Group}\lb{pfdtg}

We proceed now to project $\K^s(H_3)$ and $\K^a(H_3)$ according to the
decomposition \rep{dlrip}. Let us begin by observing that, although $\Om$ is
not trivial in $H^2(\R^2,\R/2\p)$, it is exact in another group cohomology. If
a complex of {\em gaugefied} (that is, point-dependent) k-cochains
$\R\times(\R^2)^{\pd k}\apli\R$ with a derivative $\d'$ is considered, then
there exists a 1-cochain $\Th$ such that $\Om=\d'\Th$, or
\bea
 \Om(x,y)=\d'\Th(x,y)=\Th(y\cdot q;x) - \Th(q;x+y) + \Th(q;y),
\lb{cohod}
\eea
where $y\cdot q\equiv q+y_1$ is an action of $\R^2$ on $\R$ \cite{alga}. One
such $\Th$ is given explicitly by
\bea
 \Th(q;x)=-\frac 12[(2q+x_1)x_2],
\lb{exth}
\eea
and satisfies $\Th(q;0)=0$, $\Th(q;-x)=-\Th(x^{-1}\cdot q;x)$,
$\Th(q;x)=-\Th(x\cdot q;-x)$, $\fal q\in\R,\ x\in\R^2$. By direct calculation
one also obtains the interesting property
\bea
 \Th(q;x)-\Th(y\cdot q;x)=\Th(q;-x)-\Th(y^{-1}\cdot q;-x)=y_1x_2.
\lb{nipro}
\eea
This kind of 1-cochain appears naturally in representation theory. For
example, representations \rep{idirh} on $L^2(\R)$ can be written in terms of
$\Th$ as follows \cite{var}:
\bea
 [T_{\n}(x,\al)f](q)=e^{i\n\th}e^{i\n\Th(x^{-1}\cdot
 q;x)}\,f(x^{-1}\cdot q)\hs 1\n\in\Z-\{0\}.
\lb{irrhg}
\eea

With these remarks in mind, we reinterpret the decomposition formula
\rep{dlrip}, and consider that the central element $\al=e^{i\th}\in\T$
appears in the projection in the form
\bea
 L(x,\al)\map e^{i\Th(q;-x)}\Lo(x).
\lb{prfo}
\eea
This means that we will take $\th=\Th$ and regard $\Th=\Th(q;-x)$ as a
gaugefied 1-cochain as defined above. This means that it depends on the
projected point $x\in\R^2$ which is its partner in the $H_3$ coordinates and,
furthermore, it is gaugefied -- it depends on the point $q\in\R$ where the
irreducible representations \rep{irrhg} (irreducible components of
$L(x,\al)$) act. Despite the local character of $\Th$ as regards the first
two slots $(x)$ of the $H_3$ coordinates, in the projection formula
\rep{prfo} $\Th$ will as a whole account for the third slot irrespective of
the details in its content. For example, $L(x,\al\b)$ and $L(x,\al^{-1})$
will also be projected into the right hand side of \rep{prfo}, but
$L(-x,\al)$ will be projected to $e^{i\Th(q;x)}\Lo(-x)$.

\subsubsection{The Projection of $\K^s(H_3)$}\lb{prksh3}

Let us project $\K^s(H_3)$ according to the map \rep{prfo}, in order to find
the structure of the space generated by the operators $\Lo(x),\ x\in\R^2$.
This will be done in two steps: (i) the projection of operations like norm,
involution, product, etc, and (ii) the verification of the Kac algebra axioms
for these projected operations. Here the Kac algebra axioms will just play
the role of guiding axioms, since the resulting algebra is not exactly a Kac
one.

Let's begin by projecting the norm and the involution. Since $e^{i\Th(q;-x)}$
is a complex number, from \rep{prfo} we have simply
\bean
\||L(x,\al)\||\map \||\Lo(x)\||=sup\{\||\Lo(x)\ps\||_2,\ \||\ps\||_2=1\mbox{ in
}L^2(\R^2)\}.
\eean
The dagger in $\cali M(H_3)$ is projected to
\bea
L^{\dagger}(x,\al)\map\left(e^{i\Th(q;-x)}\right)^*\,\Lo^{\dagger}(x),
\lb{invpr}
\eea
where, since these representations are unitary, $\Lo^{\dagger}(x)=\Lo(-x)$. The
involution $^*$ on the phase factor is not simply complex conjugation, but
involves also the action of $\R^2$ on $\R$. It becomes fixed if we recall that
$L^{\dagger}(x,\al)=L(-x,\al^{-1})\map e^{i\Th(q;x)}\Lo(-x)$, and compare with
\rep{invpr}, which gives
\bea
 \left(e^{i\Th(q;-x)}\right)^*=e^{i\Th(q;x)}.
\lb{invoth}
\eea
It is easy to verify that the projected $\||\ \||$ and $\dagger$ satisfy all
the usual norm and involution axioms (please, see them in Ref.~\cite{bra}).

The product in $\K^s(H_3)$ will also be projected according to \rep{prfo}. Care
must be taken when dealing with such products of operators. Since an operator
at the right {\em feels} the action of that at the left on $\R$,
its phase factor turns out to be modified. From \rep{kspr} and the above we
have
\bea
 e^{i\Th(q;-x)}\Lo(x)\,e^{i\Th(x^{-1}\cdot q;-y)}\Lo(y)=e^{i\Th(q;-x-y)}
\Lo(x+y),
\lb{prcopr}
\eea
which only gives \rep{prirp} [by \rep{cohod} and $\Om(-y,-x)=-\Om(x,y)$] if
$\Om=\d'\Th$. This implies that $\Th$ can be given by \rep{exth}. The first
condition imposed to the product is associativity, which is satisfied due to
the closeness of the cocycle $\Om$. The second, $\Lo(x)\1=\1\Lo(x)\,\fal x$,
where $\1=\Lo(0)$, is true for \rep{prirp} because $\Om(\cdot,0)=0$.

Up to this point we have a unital, involutive and normed algebra with product
given by \rep{prirp}. It is also a subalgebra of $\cali B(L^2(\R^2))$, and is
certainly closed in the weak topology defined on it. This can be seen by
comparison with the von Neumann algebra generated by the left-regular operators
$L(x)$ of $\R^2$. The only difference between the action of $\Lo(x)$ and the
action of $L(x)$ on $L^2(\R^2)$ is a phase factor (see \rep{aclo}), which does
not affect the closeness property in the weak topology (see
section~\ref{fdhg}), for example. The conclusion is that the algebra generated
by $\Lo(x),\ x\in\R^2$, is a von Neumann algebra. It will be denoted $\cali
M^{\Om}(\R^2)$.

Going further, by \rep{prfo} we project the coproduct to
\bean 
 \hat\cop L(x,\al)\map e^{i\Th(q;-x)}\copo\Lo(x),
\eean
where $\hat\cop$ is not supposed to act on the central element $e^{i\Th}$.
Taking this in the expression for the coproduct of $L(x,\al)$ and
considering again the projection formula, we get
\bea
 \copo\Lo(x)=e^{i\Th(q;-x)}\,\Lo(x)\pd \Lo(x).
\lb{coppr}
\eea
Notice that it remains symmetric, that is, $\s\circ\cop_{\Om}=\cop_{\Om}$,
where $\s(T\pd T')=T'\pd T$, as was $\hat\cop$ on $\K^s(H_3)$. The first
axiom the coproduct must satisfy is $\copo\1=\1\pd\1$, which is trivial, since
$\Th(q;0)=0$ for all $q\in\R$. The next is {\em co-associativity}, which means
\bean
 (\copo\pd id)\circ\copo=(id\pd\copo)\circ\copo.
\eean
It is also trivial since the same phase factor occurs twice in both sides of
this equation when it is applied to $\Lo(x)$. Finally, $\copo$ should be a
homomorphism from $\cali M^{\Om}(\R^2)$ to $\cali M^{\Om}(\R^2)\pd\cali
M^{\Om}(\R^2)$, which means that
\bea
 \copo(\Lo(x)\Lo(y))&=&\copo(\Lo(x))\copo(\Lo(y)).
\lb{cohom1}
\eea
The left-hand side of Equation~\rep{cohom1} yields 
\bean
 e^{i[\Om(x,y)+\Th(q;-x-y)]}\,\Lo(x+y)\pd\Lo(x+y),
\eean
while its right-hand side gives
\bean
 e^{i[\Th(q;-x)+\Th(x^{-1}\cdot q;-y)+2\Om(x,y)]}\,\Lo(x+y)\pd\Lo(x+y).
\eean
The phase factors are seen to be equal if we recall the expression for
$\Om(-y,-x)$ from \rep{cohod}, and the properties of $\Om$.

The coinvolution is projected to
\bean
 \hat\k(L(x,\al))\map e^{i\Th(q;-x)}\,\ko(\Lo(x)),
\eean
where also $\hat\k$ is supposed not to act on the phase factor. From
\rep{kscoi}, using \rep{invpr} and \rep{invoth}, we get
\bea
 \ko(\Lo(x))=e^{i[\Th(q;x)-\Th(q;-x)]} \Lo^{\dagger}(x).
\lb{coipr}
\eea
Of all the axioms imposed on a coinvolution, $\ko$ fails to satisfy only
one, the {\em anti-automorphism} axiom
\bea
 \ko(\Lo(x)\Lo(y))=\ko(\Lo(y))\ko(\Lo(x)). 
\lb{antax}
\eea
By \rep{coipr}, and taking care of the phase factors in the operator products,
we obtain from the left-hand side
\bea
 \ko(\Lo(x)\Lo(y))=e^{i[\Om(x,y)+\Th(q;x+y)-\Th(q;-x-y)]}\Lo(-x-y),
\lb{lhsant}
\eea
while the right-hand side gives
\bea
 \ko(\Lo(y))\ko(\Lo(x))=e^{i[\Th(q;y)-\Th(q;-y)+\Th(y\cdot q;x)-\Th(y\cdot
q;-x)+\Om(-y,-x)]}\Lo(-x-y).
\lb{rhsant}
\eea
After using the explicit expressions of $\Om$ and $\Th$, we get
\bea
 \ko(\Lo(x)\Lo(y))=e^{i(x_1y_2+y_1x_2)}\,\ko(\Lo(y))\ko(\Lo(x))
\lb{faiaiso}
\eea
instead of \rep{antax}. We will return to this problem below. Concerning the
remaining axioms that $\ko$ must satisfy: first, it should be involutive:
$\ko(\Lo^{\dagger}(x))=\ko(\Lo(x))^{\dagger}$. This follows from \rep{coipr}
and \rep{invpr}. The requirement $\ko(\ko(\Lo(x)))=\Lo(x)$ just
implies that the phase factor in \rep{coipr}, which is antisymmetric in $x$,
is canceled out when the second coinvolution is applied to $\Lo(-x)$. This is
obvious. The {\em anti-coautomorphism} axiom,
\bea
 \copo\circ\ko=\s\circ(\ko\pd\ko)\circ\copo,
\lb{acoiax}
\eea
is clearly satisfied: when applied to $\Lo(x)$ the left-hand side of this
equation raises the phase factor $e^{i[(\Th(q;x)-\Th(q;-x))+\Th(q;x)]}$, while
the right-hand side raises $e^{i[\Th(q;-x)+2(\Th(q;x)-\Th(q;-x))]}$, 
which is the same.

From \rep{faiaiso} it is evident that the projection $\ko$ of the
coinvolution $\k$ is not a coinvolution on $\cali M^{\Om}(\R^2)$. The role of
a coinvolution in Kac duality is explicit in formula \rep{ddi}, the
definition of the dual involution $^o$ on the predual of $\cali M(H_3)$.
Since our main goal is to prove a duality for $\cali M^{\Om}(\R^2)$, we are
faced to a serious problem. The only weak aspect of the projection process of
$\k$, which could eventually be modified to solve this problem, is the
assumption that it does not act on the phase factor. But if it did act, the
only plausible action would be by conjugation \rep{invoth} (since it acts by
{\em dagger} on $\Lo$), and the resulting $\ko$ would be just
$\ko(\Lo(x))=\Lo^{\dagger}(x)$, the usual coinvolution of a symmetric group
Kac algebra. In that case, it would not only fail to satisfy \rep{antax} but
also the anti-coautomorphism axiom, which involves the non-trivial $\copo$.
More generally, if we define $\ko$ with any phase factor other than that of
\rep{coipr}, say $e^{i\Ps(q;x)}$, the unique $\Ps$ satisfying the last three
axioms is just that combination of $\Th$'s given in \rep{coipr}. This is most
evident for the last axiom. We actually do not know of any good definition
of $\ko$ making of it a coinvolution, that is, enforcing all the above
axioms. The solution we have found for this problem is to maintain the
definition of $\ko$ as it is given by the projection, and modify the
anti-automorphism axiom \rep{antax}. A natural modification of it comes from
the projection of the anti-automorphism axiom satisfied by $L(x,\al)$, which
is given by
\bea
 \hat\k(L(x,\al)L(y,\b))=\hat\k(L(y,\b))\hat\k(L(x,\al)).
\lb{anauk}
\eea
Using the projection formula \rep{prfo} on it, we get, for example,
\bean
 \hat\k(L(x,\al)L(y,\b))\map e^{i[\Th(q;-x)+\Th(x^{-1}\cdot
 q;-y)]}\ko(\Lo(x)\Lo(y))
\eean
on its left-hand side. Doing the same with the other side, \rep{anauk} turns
out to be projected into
\bea
 \ko(\Lo(x)\Lo(y))= e^{i[\Th(q;-y)-\Th(q;-x)+\Th(y\cdot q;-x)-
\Th(x^{-1}\cdot q;-y)]}\ko(\Lo(y))\ko(\Lo(x)).
\lb{nantax}
\eea
Given its nature, \rep{nantax} should be called {\em projective
anti-automorphism} axiom. The importance of \rep{nantax} comes from the fact
that it is promptly satisfied by \rep{coipr}. In fact, substituting
\rep{lhsant} and \rep{rhsant} in \rep{nantax}, the phase factors are easily
matched with the help of the expression \rep{cohod} for $\Om$ and of its
properties. Notice that no new axiom arises if the other axioms defining a
coinvolution are projected. This ends the list of axioms satisfied by what can
now be called {\em projective coinvolution} $\ko$. It can be anticipated that
the change from the axiom \rep{antax} to \rep{nantax} will have consequences
on the predual of $\cali M^{\Om}(\R^2)$. The dual axiom, \rep{acoiax}, for the
dual coinvolution will be changed too.

Finally, the trace $\hat\ph$ is simply projected to its restriction to $\R^2$
according to
\bea
 \hat\ph(\hat f)=f(0,1)\map\pho(\hat f)=f(0),
\lb{trpr}
\eea
where a general element $\hat f$ of $\cali M^{\Om}(\R^2)$ is written 
\bea
 \hat f=\int_{\R^2} dx\ f(x)\Lo(x).
\lb{gopw}
\eea
From its very definition, this trace is n.f.s. (see section~\ref{skah3}). It
also satisfies the three specific axioms for a Haar weight, which are:
\bse
(id\pd\pho)\copo (\hat f)&=&\pho(\hat f)\1\qquad \fal \hat f\in\cali
M^{\Om}(\R^2)^+;\lb{axwe1}\slb{hwlia}\\ (id\pd\pho)[(\1\pd\hat
g^{\dagger})\copo(\hat f)]&=&\ko\circ(id\pd\pho)[\copo(\hat
g^{\dagger})(\1\pd\hat f)];\slb{hwsa}\\
\ko\circ\s^{\pho}_t&=&\s^{\pho}_{-t}\circ \ko\qquad \fal\ t\in\R.\slb{hwlic}
\ese
The third one is trivial here, since the modular group $\s^{\pho}$ is reduced
to the identity when the Haar weight $\pho$ is a trace. Concerning the other
axioms, we start by observing that, although $\pho$ is not defined on the
generators $\Lo(x)$, from \rep{trpr} and \rep{gopw} it may be guessed that it
would act on $\Lo(x)$ as $\pho(\Lo(x))=\d(x)$. This corresponds to an
extension of the domain of $\pho$ to the generators, which can be regarded as
being given by $\Lo(x)=\int_{\R^2}dy\,\d_x(y)\Lo(y)$. We notice also that
\bea
 \copo\hat f=\int_{\R^2}dx\,e^{i\Th(q;-x)}\,f(x)\,\Lo(x)\pd\Lo(x).
\lb{copfom}
\eea
Axiom \rep{hwlia} follows then from the considerations above, which imply
$(id\pd\pho)\copo\hat f=\linebreak f(0)\Lo(0)=\pho(\hat f)\1$. In the same
way the second follows from \rep{coipr} and
\bea
 \hat f^{\dagger}=\int_{\R^2}dx\,\ovl{f(x)}\,\Lo^{\dagger}(x).
\lb{hfdin}
\eea

The resulting algebra, the projection of the symmetric Kac algebra of the
Heisenberg group, will be denoted $\K^{\Om}(\R^2)$ and called {\em projective
symmetric Kac algebra} of $\R^2$. It is built on the von Neumann
algebra $\cali M^{\Om}(\R^2)$ with the usual operator norm and conjugation.
The remaining structures are grouped into:
\bse
 \Lo(x)\Lo(y)&=&e^{i\Om(x,y)}\,\Lo(x+y);\slb{kompr}\\
 \1&=&\Lo(e);\\
 \copo\Lo(x)&=&e^{i\Th(q;-x)}\,\Lo(x)\pd\Lo(x);\\
 \ko(\Lo(x))&=&e^{i[\Th(q;x)-\Th(q;-x)]}\,\Lo^{\dagger}(x);\slb{kocoi}\\
 \pho(T)&=&\left\{
\begin{array}{cl}
 \||f\||^2_2 & \mbox{if}\ T=\hat f^{\dagger}\cdot\hat f\\
 +\infty & \mbox{otherwise}
\end{array}
\right.\hs 1 T\in\cali M^{\Om}(\R^2)^+.
\ese

{\bf Comments:}
\bit
\item The product of two elements $\hat f$, $\hat g$ is
\bea
\hat f\cdot\hat g&=&\int_{\R^2\times \R^2}dx dy\
f(x)g(y)e^{i\Om(x,y)}\Lo(x+y)\non\\ &=&\int_{\R^2} dz \int_{\R^2} dx\
e^{i\Om(x,z)}f(x)g(z-x) \Lo(z)\non\\ &=&\int_{\R^2} dz\ (f\tc g)(z)\Lo(z),
\lb{poppr}
\eea
where we have used the fact that $\Om$ is antisymmetric to identify the
twisted convolution
\bea
 (f\tc g)(z)=\int_{\R^2} dx\ e^{i\Om(x,z)}f(x)g(z-x).
\lb{tcorec}
\eea 
Since the operator product is mapped into the twisted convolution of
$L^1$-functions, \rep{gopw} can be regarded as the linear left-regular
representation of $L^1_{\Om}(\R^2)$ induced by $\Lo$. $L^1_{\Om}(\R^2)$ is
the Banach algebra analogous to $L^1(\R^2)$, but with the twisted
convolution. The involution remains in $L^1$ and its image by $\Lo$ gives
the dagger of $\hat f$ (see \rep{hfdin}): $\hat f^{\dagger}=\Lo(f^*)$;
\item despite the Abelian character of $\R^2$, the projective product makes of
$\K^{\Om}(\R^2)$ a noncommutative algebra. The noncommutativity can be
measured by the introduction of a Lie algebra structure on $\cali
M^{\Om}(\R^2)$ through the commutator
\bea
 [\Lo(x),\Lo(y)]=2i\sin[\Om(x,y)]\,\Lo(x+y),
\lb{liestr}
\eea
or by the introduction of the continuous $R$ matrix 
\bean
 \Lo(x)\Lo(y)=\int_{\R^2\times\R^2} dzdw\, R(x,y;z,w)\,\Lo(z)\Lo(w).
\eean 
The $R$ matrix elements belong to $M^1_{\Om}(\R^2)\supset L^1_{\Om}(\R^2)$ and
are given by
\bea
 R(x,y;z,w)=e^{i[\Om(x,y)-\Om(z,w)]}\,\d(x+y-w-z). 
\lb{ybaesol}
\eea
Associativity of \rep{kompr} implies that $R$ should satisfy the Yang-Baxter
equation, and \rep{ybaesol} provides a new nontrivial solution for it. Let us
observe that the associativity of \rep{kompr}, the Jacobi identity for the
commutator \rep{liestr} and the Yang-Baxter equation depend on the closeness
of $\Om$;
\item the ideal of elements such that $\pho(\hat f^{\dagger}\hat f)<\infty$ is
just $\cali N_{\pho}=L^1_{\Om}\cap L^2(\R^2)$. So, the GNS representation of
this projective Kac algebra is given on $L^2(\R^2)$ by $\p_{\pho}(\hat
f)g=f\tc g$. The GNS-image in $L^2$ of $\hat f$ will be denoted $\hat
f_{\pho}$;
\item we can recover the $L^1_{\Om}$-function $f$ in the linear combination
\rep{gopw} through the Haar trace $\pho$ by the formula
\bean
 f(x)=\pho[\Lo^{\dagger}(x)\hat f];
\eean 
\item written in terms of the $M^1_{\Om}$-distributions ($L^1_{\Om}\subset
M^1_{\Om}$), the projective Kac algebra operations, other than the twisted
convolution and the trace, read ($f\in L^1_{\Om}$)
\bse
 \copo(f)(x,y)&=& e^{i\Th(q;-x)}f(x)\,\d(x-y),\\
 \ko(f)(x)&=&e^{-i[\Th(q;x)-\Th(q;-x)]}\,f(-x).
\lb{estkel}
\ese
\eit

\subsubsection{The Dual Algebra of $\K^{\Om}(\R^2)$}\lb{dkaom}

Analogously to what has been done for the Heisenberg group, we now start
fixing the predual of $\cali M^{\Om}(\R^2)$ and, by its Fourier
representation, the dual of $\K^{\Om}(\R^2)$. Even though the latter is not a
Kac algebra, the same techniques for establishing a duality for it will be
used. In this section, $L^p$ will denote $L^p(\R^2)$, for $p=1,2,\infty$. 

The elements of ${\cal M}^{\Om}(\R^2)_*$ will be written as linear forms
$\omo_{fg}$ on $\cali B(L^2)$, whose duality pairing with $\Lo^{\dagger}(x)$
gives the functions
\bean
 \bra \Lo^{\dagger}(x),\omo_{fg}\ket=(\Lo(-x)f|g)_{L^2}=(f\tc\check g)(x),\hs 1
f,g\in L^2.
\eean
The functions $\omo_{fg}(x)\equiv(f\tc\check g)(x)$ are thereby defined as the
representative functions in the predual. From this definition comes also that
these functions are in $L^{\infty}$:
$|\omo_{fg}(x)|\le\||f\||_2\||g\||_2<\infty$. When no confusion can arise, they
will be denoted simply by $f,g,h$, etc.

The product on ${\cal M}^{\Om}(\R^2)_*$ will follow from duality by
\bea
(\omo_{fg}\star\omo_{hl})(x)&=&\bra
\copo\Lo(-x),\omo_{fg}\pd\omo_{hl}\ket\non\\
&=&e^{i\Th(q;x)}\,(\Lo(-x)\pd\Lo(-x)(f\pd h)|g\pd l)_{L^2\pd L^2}\non\\
&=&e^{i\Th(q;x)}\,\omo_{fg}(x)\omo_{hl}(x),
\lb{stpr}
\eea
and, since $\copo$ is symmetric, $\star$ is Abelian. Its associativity follows
trivially. The unit is a consequence of \rep{stpr} and is uniquely given by
\bean
 \1(x)=e^{-i\Th(q;x)}.
\eean

To finish with the characterization of this predual, the dual involution
$^o$ is find out through
\bea
(\om_{fg}^{\Om})^o(x)&=&\bra\Lo(-x),\om^{\Om\ o}_{fg}\ket=\ovl{\bra\ko(\Lo^{\dagger}(-x)),\omo_{fg}\ket}\non\\
 &=&e^{-i[\Th(q;x)-\Th(q;-x)]}\ovl{\omo_{fg}(x)}.
\lb{predinv}
\eea
It is indeed an involution, since it is antilinear, satisfies
$^o\,\circ\,^o=id.$, $\1^o(x)=\1(x)\ \fal x$, and is an anti-automorphism
because, after recalling \rep{invoth}, the two lines below turn out to be
equal:
\bean
(f\star g)^o(x)&=&e^{-i[\Th(q;x)-\Th(q;-x)]}\left(e^{i\Th(q;x)}\right)^*
\ovl{f(x)g(x)}\\
(g^o\star f^o)(x)&=&e^{i\Th(q;x)-2[\Th(q;x)-\Th(q;-x)]}\ovl{g(x)f(x)}.
\eean

In analogy with the Kac algebra case, we call the predual ${\cal M}^{\Om}(\R^2)_*$ the {\em projective Fourier algebra} of $\R^2$ and denote it by
$A^{\Th}(\R^2)$. Since their elements are $L^{\infty}$, we will denote the
respective von Neumann algebra by $L^{\infty}_{\Th}(\R^2)$. As a final remark,
notice that the projection of the product $(f\cdot
g)(x,\al)=f(x,\al)g(x,\al)$, giving the correct projected product \rep{stpr},
comes from the projection of $f\in L^{\infty}(H_3)$ to
$L^{\infty}_{\Th}(\R^2)$ in the following way:
\bea
 f(x,\al)\map e^{-i\Th(q;-x)}f(x).
\lb{prlif}
\eea
Actually, it gives the $\star$ product if we allow the phase $\Th$ to {\em
feel} the action of $\R^2$ on $\R$ as follows: $f(x,\al)g(y,\b)=(f\pd
g)(x,\al\pd y,\b)$ is projected, after \rep{prlif} and considering that $y$
acts on the phase at $x$, to $e^{-i\Th(y\cdot
q;-x)}e^{-i\Th(q;-y)}\,f(x)g(y)$, while $(f\cdot g)(x,\al)$ goes to
$e^{-i\Th(q;-x)}\,(f\star g)(x)$. Making $(y,\b)=(x,\al)$ and recalling that
$-\Th(x\cdot q;-x)=\Th(q;x)$, the exponentials without action cancel out and
the correct $\star$ product is obtained. In the same way, the action of
$L^{\infty}(H_3)$ on the Hilbert space $H(H_3)$ by pointwise product is
projected to the action of $L^{\infty}_{\Th}(\R^2)$ on $L^2(\R^2)$ by the
$\star$ product. Projection \rep{prlif} is also useful to project
$L^1$-functions. From \rep{prlif} and \rep{prfo}, we get the right projection
of the $\K^s(H_3)$-elements \rep{lrrl} into the operators \rep{gopw}.

The Fourier representation $\l^{\Om}$ of $A^{\Th}(\R^2)$ is defined, in analogy
with \rep{frdef}, by
\bea
 [\l^{\Om}(\omo)f](x)=[(\omo\circ\ko\pd id)\copo\hat f]_{\pho}(x)\hs 1 f=\hat
f_{\pho}\in L^2\cap L^1_{\Om}.
\lb{prdeffr}
\eea
From \rep{copfom} and \rep{kocoi} we readily get
\bean
(\omo\circ\ko\pd id)\copo\hat f&=&(\omo\pd
id)\int_{\R^2}dx\,e^{i\Th(q;x)}f(x)\,\Lo(-x)\pd\Lo(x)\\ &=&
\int_{\R^2}dx\,e^{i\Th(q;x)}f(x)\omo(x)\,\Lo(x),
\eean
which gives 
\bean
 \l^{\Om}(\omo)f(x)=(\omo\star f)(x).
\eean
Since the action of $A^{\Th}(\R^2)$ (or $L^{\infty}_{\Th}(\R^2)$) is given by
the $\star$ product, it follows that $\l^{\Om}=id$. This means that, as in
the Kac case, the dual of $\K^{\Om}(\R^2)$ is built on the von Neumann
algebra $L^{\infty}_{\Th}(\R^2)$, the $L^{\infty}$-Banach algebra with the
product $\star$ and the involution $^o$ of $A^{\Th}(\R^2)$. Its norm is the
same of that $L^{\infty}$ and satisfies $\||f^o\||=\||f\||$, $\||f^o\star
f\||=\||f\||^2$.

The generator of $\l^{\Om}$, that is, the operator $W^{\Th}$ in
$L^{\infty}_{\Th}\pd\cali B(L^2)$ satisfying
$\l^{\Om}(\omo)=(id\pd\omo)(W^{\Th})$, is given by
\bea
 [W^{\Th}F](x,y)=e^{i\Th(q;x)}e^{-i\Om(x,y)}\,F(x,y+x).
\lb{wth}
\eea
We see that this generator acts on $(x,y)\in\R^2\times\R^2$ as
if it were $W^{\Th}\sim 1\pd\Lo^{\dagger}(x)$, with $1$ the constant function
$1\in L^{\infty}_{\Th}$.

By duality, we obtain also a coproduct on the dual $L^{\infty}_{\Th}(\R^2)$:
\bea
\copt(\omo_{fg})(x,y)&=&\bra [\Lo(x)\Lo(y)]^{\dagger},\omo_{fg}\ket
=e^{-i\Om(x,y)}(\Lo^{\dagger}(x+y)f|g)\non\\
&=& e^{-i\Om(x,y)}\omo_{fg}(x+y).
\lb{prcopt} 
\eea
This coproduct is automatically coassociative, as is defined by the
associative dual product. It is also unital:
$\copt\1(x,y)=e^{-i\Om(x,y)}\,e^{-i\Th(q;x+y)}$, while
$(\1\pd\1)(x,y)=e^{-i[\Th(y\cdot q;x)+\Th(q;y)]}$. In this last expression it
must be recalled that the phase factor $e^{i\Th}$ is a special kind of
complex function which is {\em sensitive} to the action of $\R^2$ on $\R$,
even when the simple product in $\C$ is performed. The homomorphism axiom was
already proved for $\K^{\Om}(\R^2)$, but it is interesting to verify it here
again, since it confirms the strange behavior of $e^{i\Th}$ under the complex
product. It follows from:
\bean
 \copt(f\star g)(x,y)&=&e^{-i\Om(x,y)}\,e^{i\Th(q;x+y)}\,f(x+y)g(x+y)\\
 (\copt f\star\copt g)(x,y)&=&(\copt f^{(1)}\star\copt g^{(1)})(x)(\copt
f^{(2)}\star\copt g^{(2)})(y)\\ &=&e^{i\Th(y\cdot q;x)}\,e^{i\Th(q;y)}\,\copt
f(x,y)\copt g(x,y)\\ &=&e^{i\Th(y\cdot
q;x)}\,e^{i\Th(q;y)}\,e^{-2i\Om(x,y)}\,f(x+y)g(x+y),
\eean
where we have written $\copt f=\copt f^{(1)}\pd\copt f^{(2)}$. The
homomorphism is established after recalling the expression for $-\Om(x,y)$.
The above two axioms would not hold if we did not allow $\Th$ to {\em feel} 
the action of $\R^2$ on $\R$.

The candidate coinvolution in $L^{\infty}_{\Th}$ comes from the
duality relation
\bea
 \kt(\omo_{fg})(x)=\bra\ko(\Lo(-x)),\omo_{fg}\ket=
 e^{-i[\Th(q;x)-\Th(q;-x)]}\omo_{fg}(-x).
\lb{ducoi}
\eea
It satisfies the first three coinvolution axioms without any problem,
including the anti-isomor\-phism condition which its dual $\ko$ fails to
satisfy. The problem lies precisely in the anti-coauto\-mor\-phism axiom.
This is just the dual of the problem found in $\K^{\Om}(\R^2)$, and its
solution will be given by dualizing the solution of that problem, namely
dualizing the axiom Eq.~\rep{nantax}. Recall that \rep{antax} can be written
in the form $\ko\circ m=m\circ(\ko\pd\ko)\circ\s$, where $m$ denotes the
product on $\K^{\Om}(\R^2)$. We adapt this form to axiom \rep{nantax} and
dualize it, that is, just transpose the order and change the operations to
their duals, taking into account the effect ~\rep{invoth} of duality on $\Th$
and transposing $x\leftrightarrow y$. It should be noticed that, from the
properties of $\Th$, it follows also that $(e^{i\Th(x^{-1}\cdot
q;y)})^*=e^{i\Th(x\cdot q;-y)}$. We obtain
\bea
 [\copt\circ\kt f](x,y)= e^{-i[\Th(q;y)-\Th(q;x)+\Th(y\cdot q;x)-
\Th(x^{-1}\cdot q;y)]}\, [\s\circ(\kt\pd\kt)\circ\copt f](x,y)
\lb{pracoax}
\eea
as a new axiom replacing \rep{acoiax}. It will be called {\em projective
anti-coautomorphism} axiom. Provided care is taken with the product of 
complex functions $e^{i\Th}$ when computing its right-hand side, \rep{ducoi} is
promptly seen to satisfy this new axiom. After these changes, $\kt$ should be
called {\em coprojective coinvolution}.

Finally, the trace \rep{kaph} is projected to the n.f.s. trace
\bea
\ph(f)=\int_{H_3}dxd\al\,f(x,\al)\map\pht(f)=\int_{\R^2}dx\,e^{-i\Th(q;-x)}
f(x),
\lb{prtrt}
\eea
where $f\in L^{\infty}_{\Th}(\R^2)^+$. If this projection is interpreted as
coming from the projection \rep{prlif} of $f\in L^{\infty}(H_3)$ to
$L^{\infty}_{\Th}(\R^2)$, it is evident that the projection of $\ph$ to
$\pht$ is only possible due to the compactness of the central group $\T$. This
fact has also been observed in prequantization \cite{tuy}.

The trace \rep{prtrt} is of Haar type, since it satisfies the axiom \rep{hwsa},
for example, as follows: using the formulas for $\pht,\ ^o$ and $\star$, the 
left-hand side gives
\bean
 (id\pd\pht)[(\1\pd g^o)\star\copt
f](x)&=&\int_{\R^2}dy\,e^{-i\Th(q;-y)}\,e^{i\Th(q;y)}g^o(y)\copt f(x,y)\\
&=&\int_{\R^2}dy\,\copt f(x,y)\ovl{g(y)},
\eean
while from the right-hand side we have
\bean
 \kt\{(id\pd\pht)[\copt g^o\star(\1\pd f)]\}(x)&=&
e^{-i[\Th(q;x)-\Th(q;-x)]}(\copt g^o)^{(1)}(-x)\pht[(\copt g^o)^{(2)}\star f]\\
&=&\int_{\R^2}dz\, \ovl{\copt g(-x,z)}f(z).
\eean
They become equal when we substitute the expression for the coproduct $\copt$
and make the change of variable $y=z-x$ in the last integral. Axiom
\rep{hwlia} follows similarly and \rep{hwlic} comes automatically because
$\s^{\pho}$ reduces to the identity.

The algebra so far obtained will be called the {\em projective Abelian Kac
algebra} of $\R^2$ and will be denoted $\K^{\Th}(\R^2)$. It is built on
$L^{\infty}_{\Th}(\R^2)$ and its structure, besides the usual
$L^{\infty}$-norm, is summed up in the following properties:
\bse
 f^o(x)&=&e^{-i[\Th(q;x)-\Th(q;-x)]}\,\ovl{f(x)}\slb{invprab}\\
 (f\star g)(x)&=&e^{i\Th(q;x)}\,f(x)g(x)\slb{abstpr}\\
 \1(x)&=&e^{-i\Th(q;x)}\\
 \copt f(x,y)&=&e^{-i\Om(x,y)}\,f(x+y)\\
 \kt f(x)&=&e^{-i[\Th(q;x)-\Th(q;-x)]}\,f(-x)\\
 \pht(f)&=&\int_{R^2}dx\, e^{-i\Th(q;-x)}\,f(x).
\ese

{\bf Comments:}
\bit
\item the name {\em projective} does not put $\K^{\Om}(\R^2)$ and
$\K^{\Th}(\R^2)$ into the same category, since the {\em projective}
coinvolution $\ko$ and the {\em coprojective} coinvolution $\kt$ are defined
by different sets of axioms;
\item it can be confirmed that $L^{\infty}_{\Th}(\R^2)$ acts on $L^2(\R^2)$
from the GNS construction induced by the trace $\pht$. This trace defines a
scalar product on the ideal of elements $f$ such that $\pht(f^o\star f)<\infty$
by $(f|g)\equiv\pht(g^o\star f)$. Direct calculation from \rep{invprab} and
\rep{abstpr} gives the usual $L^2$-scalar product $(f|g) =\int_{\R^2}dx,
f(x)\ovl{g(x)}$. The inclusion of $f\in L^{\infty}_{\Th}$ into $L^2$ by this
construction is denoted $f_{\pht}$.
\eit

The scalar product defined by $\pht$ also shows that the involution $^o$ goes
to the simple complex conjugation in $L^2$. Thus, the antiunitary operators
implementing on $L^2$ the involutions $^o$ and $^*$ are simply
\bean
 J^{\Th}f(x)&=&\ovl{f(x)}\\
 J^{\Om}f(x)&=&\ovl{f(-x)}.
\eean 
They provide the adjoint of the unitary Fourier representation generators, for
example,
\bea
W^{\Th *}=(J^{\Om}\pd J^{\Th})W^{\Th}(J^{\Om}\pd J^{\Th}),
\lb{jconj}
\eea
whose action on $L^2\times L^2$ is 
\bean
 [W^{\Th *}F](x,y)=e^{i\Th(q;-x)}e^{i\Om(x,y)}\,F(x,y-x).
\eean

Operator $J^{\Om}$ is the correct projection of $\hat J$, but $J^{\Th}$ is
just the implementation of $^o$ as complex conjugation and does not equal
the projection of $[Jf](x,\al)=\ovl{f(x,\al)}$. The latter is projected to
\bean
 J'^{\Th}f(x)=e^{-i[\Th(q;x)-\Th(q;-x)]}\ovl{f(x)}.
\eean
In terms of the actual projected operators we have the canonical
implementations of $\ko$ and $\kt$ in $L^2$:
\bean
 \ko(\Lo(x))&=&J'^{\Th}\,\Lo^{\dagger}(x)\,J'^{\Th}\\
 \kt(f)&=&J^{\Om}\,f^o\,J^{\Om}.
\eean
That $\ko$ is implemented by $J'^{\Th}$, and not by $J^{\Th}$, can be explained
if we recall that: 1) differently from $\kt$, $\ko$ is not a linear
anti-automorphism but a projective one; 2) the involutions embodied in the
antiunitary operators $J$ are closely related to the algebra product and not
to the coproduct.

As regards the canonical implementation of $\copt$ and $\copo$, it is easily
verified that $W^{\Th}$ and $W^{\Om}$, respectively, do the job. After
recalling the $\star$-action of $L^{\infty}_{\Th}$ on $L^2$, it is found that
\bean 
[W^{\Th}(\1\pd f)W^{\Th *}F](x,y)=e^{i[\Th(q;x)+\Th(q;x+y)+\Th(x\cdot
q;-x)]}f(x+y)F(x,y).
\eean
The phase $\Th(x\cdot q;-x)$ comes from the action of $W^{\Th *}$ over the
previous action of $W^{\Th}$ by $\Lo(-x)$, and cancels out with the first
$\Th$. On the other hand, we have
\bean
 [\copt(f)F](x,y)=e^{-i[\Om(x,y)+\Th(y\cdot q;x)+\Th(q;y)]}f(x+y)F(x,y).
\eean
Recalling the expression for $\Om$ as the cohomological derivative of
$\Th$, we arrive at the equality of the left-hand sides, namely
\bea
\copt(f)=W^{\Th}(\1\pd f)W^{\Th *}.
\lb{cicopt}
\eea
This implementation is not unique, as \rep{cicopt} will be also satisfied for
any operator $V^{\Th}$ introduced through
$W^{\Th}F(x,y)=e^{i\Th(q;x)}\,V^{\Th}F(x,y)$ or, equivalently, by
$V^{\Th}F(x,y)=e^{-i\Om(x,y)}F(x,y+x)$. $W^{\Th}$ is, nevertheless, the
unique operator generating the Fourier representation $\l^{\Om}$. In Kac
algebra terminology, it is the {\em fundamental operator} of the projective
Kac algebra $\K^{\Th}(\R^2)$. This means that it satisfies the pentagonal
relation, a point easily confirmed by direct calculation.

To show the implementation of $\copo$, we first find
$W^{\Om}=\widehat{W^{\Th}}=\s\circ W^{\Th *}\circ\s$, and $W^{\Om *}=\s\circ
W^{\Th}\circ\s$,
\bse
 W^{\Om}F(x,y) &=& e^{i\Om(y,x)}e^{i\Th(q;-y)}\,F(x-y,y)\slb{wom}\\
 W^{\Om *}F(x,y) &=& e^{-i\Om(-y,x)}e^{i\Th(q;y)}\,F(x+y,y).
\ese
Keeping in mind that these generators behave like $W^{\Om}\sim\Lo(y)\pd 1^o$
and $W^{\Om *}\sim\Lo^{\dagger}(y)\pd 1$ we get, after some cancellations,
\bean
 [W^{\Om}(I\pd\Lo(z))W^{\Om *}F](x,y)=
 e^{i[\Om(z,x)+\Th(q;-y)+\Th(y^{-1}\cdot q;y-z)]}F(x-z,y-z)
\eean
and
\bean 
 \copo\Lo(z)F(x,y)=e^{i[\Th(q;-z)+\Om(z,x)+\Om(z,y)]}F(x-z,y-z),
\eean
which turn out to be equal if we recall that $\Om(z,y)=\Om(y-z,-y)$. Unlike
what happens in the dual case, the operator defined from the expression for
$W^{\Om}$ by $V^{\Om}F(x,y)=e^{-i\Om(x,y)}F(x-y,y)$ does not implement $\copo$.

In order to establish a projective duality, we proceed now to determine the
predual of $L^{\infty}_{\Th}$ and its Fourier representation. Before that, we
must point out a particularity of projective algebras: the duality pairing
$\bra ,\ket$, used so far to connect $\cali M^{\Om}(\R^2)$ and its predual
$A^{\Th}(\R^2)$, involves implicitly a complex conjugation of the phase
factors, since the representative functions in the projective Fourier algebra
are defined in \rep{dplf} by pairing with $\Lo^{\dagger}$ and not with $\Lo$.
The result is that the phase factors in the structure of $A^{\Th}(\R^2)$ get
complex-conjugated. Furthermore, the dualization process from
Eq.~\rep{pracoax}, which brings the projective anti-automorphism axiom from
$\K^{\Om}(\R^2)$ to the projective anti-coautomorphism axiom of
$\K^{\Th}(\R^2)$, involves not only a transposition but also a complex
conjugation of the phase factors. Differently from these duality pairings,
the impossibility of working with the $L^{\infty}_{\Th}$-generators forces us
to put generic $L^{\infty}_{\Th}$-functions $g$ in the duality pairing
between this algebra and its predual, which is usually given by $\bra
g,f\ket=\int_{\R^2}dx\,g(x)f(x)$. That is, the complex conjugation of the
phase factors is {\em implicit}. With these remarks in mind, and recalling
that the $\Th$-phase factors are complex conjugated according to
\rep{invoth}, while those involving $\Om$ are complex conjugated as usual, we
begin by introducing the representative functions on $L^{\infty}_{\Th *}$
through
\bea
 \bra g,\omt_{hf}\ket&=&(g\star h|f)=\int
dx\,g(x)\,e^{i\Th(q;-x)}h(x)\ovl{f(x)},\hs 1 g\in L^{\infty}_{\Th},\ h,f\in
L^2\non\\ &\therefore&\omt_{hf}(x)=e^{i\Th(q;-x)}h(x)\ovl{f(x)}=(h\star
f^o)(x).
\lb{repdf}
\eea
When no confusion can arise these functions will be also denoted by $f,g,h$,
etc. The product in the predual comes from the coproduct $\copt$ by
\bean
 \bra g,f\tc h\ket=\bra\copt g, f\pd
 h\ket=\int_{\R^2\times\R^2}dx\,dy\,e^{i\Om(x,y)} g(x+y)f(x)h(y),
\eean
which gives the twisted convolution \rep{tcorec}. The involution $^*$ comes
from
\bean
 \bra g,f^*\ket=\ovl{\bra\kt(g^o),f\ket}=\int_{\R^2}dx\,g(x)\ovl{f(-x)}, 
\eean
and coincides with the $L^1$-involution. At this point we could already guess
that the predual we are looking for is just $L^1_{\Om}(\R^2)$. This is
confirmed when we recall the H\"older inequality \cite{3sen} for
$L^p$-spaces, $p=1,2$, which says that, if $h,f\in L^2$, then the modulus of
their product is an integrable function in $L^1$. This implies that the
functions $\omt_{hf}$ given in \rep{repdf} are $L^1$. Furthermore, their
product and involution also characterize them as $L^1_{\Om}(\R^2)$-functions.

If duality is to hold, the dual of the Fourier representation $\l^{\Om}$ should
be generated by $W^{\Om}$, the dual of $W^{\Th}$. This dual generator has
already been given in \rep{wom}. The representation of
$L^1_{\Om}(\R^2)$ it generates is denoted by $\l^{\Th}$ and follows from the
identity
\bean
 (g|\l^{\Th}(\omt_{hf})\star l)=(g\pd f|W^{\Om}(l\pd h)).
\eean
Recalling that $\omt_{hf}=h\tc\check f\in L^1_{\Om}$, we find that
$\l^{\Th}$ is given by
\bea
 \l^{\Th}(\omt)=\int_{\R^2}dx\,\omt(x)\,\Lo(x).
\lb{frepth}
\eea
Its range as an operator on a Hilbert space is restricted to
$L^{\infty}_{\Th}\cap L^2(\R^2)$. Let us observe that $\l^{\Th}$ cannot be
written $\l^{\Th}(\omt)f_{\pht}=[(\omt\circ\kt\pd id)(\copt)(f)]_{\pht}$,
perhaps because $\kt$ is not an anti-coautomorphism.

Formula \rep{frepth} coincides with the expression \rep{gopw} for a generic
element of $\cali M^{\Om}(\R^2)$, which enables us to conclude, from
\rep{poppr} and \rep{hfdin}, that it is actually a linear and involutive
representation of $L^1_{\Om}(\R^2)$ in that von Neumann algebra. At this
point it is no more necessary to show that the product $\star$ and the
coprojective coinvolution $\kt$ go, by duality, respectively into the
coproduct $\copo$ and the projective coinvolution $\ko$. Furthermore, in
addition to being the unique operator implementing the coproduct $\copo$ and
generating $\l^{\Th}$, the fundamental operator $W^{\Om}$ satisfies the
pentagonal relation. All these facts confirm the existence of a duality
between the projective Kac algebras $\K^{\Om}(\R^2)$ and $\K^{\Th}(\R^2)$. By
the association of these symmetric and Abelian projective Kac algebras to the
Abelian group $\R^2$, the projective Kac duality provides a {\em projective
Fourier duality} for this group.

As a by-product of the the pentagonal relation, which can be
considered as {\em the symbol of duality} \cite{ska} and is satisfied by both 
$W^{\Th}$ and $W^{\Om}$, we find that the operators $V^{\Th}$ and $V^{\Om}$, 
coming from
\bean
 && W^{\Th}F(x,y)=e^{i\Th(q;x)}\,V^{\Th}F(x,y),\\
 && W^{\Om}F(x,y)=e^{i\Th(q;-y)}\,V^{\Om}F(x,y),
\eean
satisfy the following {\em projective} versions of that relation:
\bean
&&[V^{\Th}_{23}V^{\Th}_{13}V^{\Th}_{12}F](x,y,z)=e^{-i\Om(x,y)}\,
[V^{\Th}_{12}V^{\Th}_{23}F](x,y,z),\\
&&[V^{\Om}_{23}V^{\Om}_{13}V^{\Om}_{12}F](x,y,z)=e^{-i\Om(y,z)}\,
[V^{\Om}_{12}V^{\Om}_{23}F](x,y,z).
\eean

As regards projective Kac duality, it must be remembered that these algebras
are not objects in the same category if the Kac algebra category definition
given in Ref.~\cite{ensc} is to be maintained. The coinvolutions in
the projective Kac algebras satisfy different sets of axioms. They would
become objects of the same category if we could define a wider category whose
objects would be algebras similar to Kac algebras, but where the coinvolutions 
would be more general linear maps $\k'$ satisfying only the axioms
\bean
 \k'\circ \ast &=&\ast\circ\k'\\
 \k'\circ\k' &=&id.
\eean
Unfortunately, such a category is not well defined, since the
anti-coauto\-mor\-phism property of $\ko$ and the anti-automorphism property of
$\kt$ seem to play an important role in projective Kac duality. For example,
the alluded property of $\ko$ seems to be the responsible for the expression
\rep{prdeffr}, while the same is not true between $\kt$ and $\copt$.

\subsection{Irreducible Decomposition According to the Projective
Dual}\lb{idska}

This subsection is devoted to the decomposition of the projective Kac duality
obtained in the last subsection according to the projective unitary dual of
$\R^2$. We begin with the decomposition of the projective Kac algebra
$\K^{\Om}(\R^2)$ according to the decomposition of the left-regular
representations $\Lo$ in terms of projective irreducible representations. The
latter are obtained by restricting to $\R^2$ the irreducible linear
representations of $H_3$ shown in
\rep{stonevn}. The result is
\bea
 [S_{\n}(x)\x](q)=e^{-i\n\Th(q;-x)}\x(q-x_1), \hs 1\n\in\Z-\{0\}.
\lb{irprorep}
\eea
By direct calculation and using the relation between the cochains $\Th$ and
$\Om$ one easily verifies that these operators satisfy the projective
relation \rep{popr}, while the members of the other series of irreducible
representations of $H_3$, \rep{irchar}, do not. By the Stone-von Neumann
theorem and Bargmann's method, we conclude that \rep{irprorep} are the unique
irreducible projective representations of the plane group. Since they are
also inequivalent, the $\Om$-projective dual of $\R^2$, here denoted
$\widehat{\R^2_{\Om}}$, is just $\Z-\{0\}$. The von Neumann algebra generated
by this kind of bounded operators on $L^2(\R)$ will be denoted $\cali
M^{\Om}_{\n}(\R^2)$. Since the operators \rep{irprorep} come from the
representations of the Heisenberg group, which is a group of type~I, the
algebra $\cali M^{\Om}_{\n}(\R^2)$ is also of type~I. In the following we
will proceed along the lines of the symmetric Kac algebra decomposition
exposed in Ref.~\cite{rasa}. The task here will be simpler than in that work,
since the Haar weight involved is a trace. The decomposition of a von Neumann
algebra generated by regular representations of a unimodular type~I group was
already established in Ref.~\cite{dix2}. The only new aspect here is that the
representations involved are {\em projective}. To proceed further, we will
suppose the existence of a positive measure $\m(\n)$ on
$\widehat{\R^2}_{\Om}$ such that the following equality holds true:
\bea
 \Lo=\sum^{\+}_{\n\in\Z-\{0\}}\m(\n)\, S_{\n}. 
\lb{prlrrde}
\eea
 The decomposition \rep{prlrrde} is based on the facts that (i) both
$\Lo$ and $S_{\n}$ are projective operators (ii) these are irreducible; (iii)
the representations
\rep{irprorep} can also be defined by
\bea
 [S_{\n}(x)f_{\n}](y)=e^{i\n\Om(x,y)}\,f_{\n}(y-x),
\lb{defprpi}
\eea
where the functions $f_{\n}\in H_{\n}(\R^2)\sim L^2(\R)$ enter in the
projective decomposition of $f\in L^2(\R^2)$ according to
$f=\sum_{\n\in\Z-\{0\}}\m(\n)\,f_{\n}$, and are given by
$f_{\n}(x)=e^{\frac{i\n}2 x_1x_2}\x(x_1)$. Taking this Hilbert-space 
projective decomposition into \rep{defprpi}, we promptly obtain
\rep{irprorep}, with $q=y_1$. Thus, $S_{\n}$ is actually an irreducible
projective component of $\Lo$.

Since \rep{prlrrde} does make sense, we can go on and decompose the
respective representation of $L^1_{\Om}(\R^2)$ through
\bean
 \Lo(f)=\sum_{\n\in\Z-\{0\}}\m(\n)\, S_{\n}(f),
\eean
where 
\bea 
 S_{\n}(f)\equiv \hat f_{\n}=\int_{\R^2} dx\, f(x)\, S_{\n}(x).
\lb{preweyl}
\eea
The last formula can be regarded as the {\em projective Fourier
transform} on $\R^2$, that is, a map associating an operator-valued
function of $\widehat{\R^2}_{\Om}$ to each $L^1_{\Om}$-function on $\R^2$.

Operators \rep{preweyl} act on $L^2(\R)$ according to \rep{irprorep}, through
\bean
 [\hat f_{\n}\x](q)=\int_{\R} du\,K^{\n}_f(q,u)\,\x(u),
\eean
where the kernel $K^{\n}_f$ is given by
\bean
 K^{\n}_f(q,u)=\int_{\R} dv\,e^{-i\n\Th(q;(u-q,-v))}\,f(q-u,v).
\eean
This enables us to introduce a trace on the operators \rep{preweyl} by
\bean
 \Tr_{\n}(\hat f_{\n})\equiv\frac 1{2\pi\n}\int_{\R} dq\, K_f^{\n}(q,q).
\eean
After recognizing the Dirac delta distribution on $\R$,
$\d(v)=\frac 1{2\pi}\int_{\R}dq\,e^{iqv}$, the trace above turns out to be
\bea
 \Tr_{\n}(\hat f_{\n})=f(0)
\lb{prttr}
\eea
for all $\n\in\Z-\{0\}$. On the other hand, by Ref.~\cite{dix2}, the
decomposition of the trace $\pho$ on $\K^{\Om}(\R^2)$ should be given by
n.f.s. traces $\pho_{\n}$ according to
\bea
 \pho(\hat f)=\sum_{\n\in\Z-\{0\}}\m(\n)\,\pho_{\n}(\hat f_{\n}).
\lb{trdec}
\eea
If we take $\pho_{\n}=\Tr_{\n}$ and recall that $\pho(\hat f)=f(0)$, we
conclude from \rep{prttr} that the measure $\m$ must be such that
$\sum_{\n\in\Z-\{0\}}\m(\n)=1$. As a first application of the trace
decomposition, we establish a projective version of the Plancherel formula.
It comes as a consequence of $\pho(\hat f^{\dagger}\hat f)=(f^*\tc f)(0)$ and
of the decomposition \rep{trdec},
\bea
\int_{\R^2}dx\, |f(x)|^2=\sum_{\n\in\Z-\{0\}}\m(\n)\,\Tr_{\n}[\hat
f_{\n}^{\dagger}\hat f_{\n}],
\lb{prplan}
\eea
which allows $\m$ to be called the {\em projective Plancherel measure}
associated to the Haar measure on $\R^2$.

We recover the $L^1_{\Om}$-function at \rep{preweyl}, or the inverse
projective Fourier transform, by
decomposing $f(x)=\pho[\Lo^{\dagger}(x)\hat f]$, namely,
\bea
 f(x)&=&\sum_{\n\in\Z-\{0\}}\m(\n)\,\Tr_{\n}[S^{\dagger}_{\n}(x)\hat
f_{\n}]\non\\ &\equiv&\sum_{\n\in\Z-\{0\}}\m(\n)\,f_{\n}(x).
\lb{l1fdec}
\eea
Since $f_{\n}(x)=\Tr_{\n}[S^{\dagger}_{\n}(x)\hat f_{\n}]=f(x)$ for all $\n$
(see \rep{prttr}), the sum over the dual is, in fact, not needed and we have
\bea
 f(x)=\Tr_{\n}[S^{\dagger}_{\n}(x)\hat f_{\n}],\hs 1 \fal\n\in\Z-\{0\}.
\lb{wigrec}
\eea

The von Neumann algebra $\cali M^{\Om}_{\n}(\R^2)$ with the projective
operator product of its generators $S_{\n}(x),\ x\in\R^2$, together with the
trace $\Tr_{\n}$ and the remaining projective Kac algebra structure inherited
from $\K^{\Om}(\R^2)$, turns out to be a projective Kac algebra. This is so
because, by property \rep{prttr}, the traces $\Tr_{\n}$ have the same
characteristics of the trace $\pho$, and are also Haar traces. This algebra
will be denoted $\K^{\Om}_{\n}(\R^2)$, and its structure is given by
\bse
 S_{\n}(x)S_{\n}(y)&=&e^{i\n\Om(x,y)}\,S_{\n}(x+y);\\
 \1&=&S_{\n}(0);\\
 \cop^{\n}S_{\n}(x)&=&e^{i\n\Th(q;-x)}\,S_{\n}(x)\pd S_{\n}(x);\slb{komncop}\\ 
 \k^{\n}S_{\n}(x)&=&e^{i\n[\Th(q;x)-\Th(q;-x)]}\,S_{\n}^{\dagger}(x);\\
 \Tr_{\n}(T_{\n})&=&\left\{
\begin{array}{cl}
 \|| f\||^2_2 & \mbox{if}\ T_{\n}=\hat f_{\n}^{\dagger}\cdot\hat f_{\n}\\
 +\infty & \mbox{otherwise}
\end{array}
\right.\hs 1 T_{\n}\in\cali M^{\Om}_{\n}(\R^2)^+.
\ese 
These projective Kac algebras have exactly the same structure as that of
$\K^{\Om}(\R^2)$, so it is unnecessary to verify the axioms again. 
Observe that its GNS-representation, as induced by $\Tr_{\n}$, is in the
Hilbert space $L^2(\R^2)$, while its elements act on $L^2(\R)$. This is due to
the fact that the operators
\bea
 \hat f_{\n}=\int_{\R^2}dx\, f(x)\,S_{\n}(x)
\lb{eldekomn}
\eea
are written in terms of $L^1_{\Om}(\R^2)$-functions $f$, while the generators
$S_{\n}$ act on the ``wavefunctions'' on the configuration space.

The main difference between the above projective decomposition and the linear
one performed in Ref.~\cite{rasa} lies in the Haar weight decomposition. In
the linear case the Haar weight (trace or not) does not satisfy \rep{prttr}
and is not, consequently, decomposed into Haar weights as happens with
$\pho$. In that case the Kac algebra decomposes into {\em Hopf-von Neumann}
algebras generated by irreducible operators, and not into algebras of the
same category (recall that a Kac algebra is just a Hopf-von Neumann algebra
plus a n.f.s.\frs Haar weight).

The predual of $\cali M_{\n}^{\Om}(\R^2)$ is obtained in the same way as
in the previous case, that is, by duality. The representative
functions $\oml_{\x\ch}$ in the predual are given by the pairing
\bea 
\oml_{\x\ch}(x)&\equiv&\bra S_{\n}^{\dagger}(x),\oml_{\x\ch}\ket=(S_
{\n}^{\dagger}(x)\x|\ch)_{L^2(\R)}\non\\
&=&\int_{\R} dq e^{-i\n\Th(q;x)}\,\x(q+x_1)\ovl{\ch(q)}.
\lb{prewig}
\eea
Also by duality, we obtain that the involution is conjugation by $^o$ and the
product is the star product $\star$, operations already introduced in
\rep{predinv} and \rep{stpr}. The only difference between the
present operations and those previously shown is a $\n$-dependence in the phase
factors. They are given explicitly by
\bean
 f_{\n}^o(x)&=&e^{-i\n[\Th(q;x)-\Th(q;-x)]}\,\ovl{f_{\n}(x)}\\
 (f_{\n}\star g_{\n})(x)&=&e^{i\n\Th(q;x)}\,f_{\n}(x)g_{\n}(x),
\eean
where we have written the functions \rep{prewig} as $f_{\n},\ g_{\n}$, etc.,
to emphasize their $\n$-dependence. From their definitions it follows also
that these functions are essentially bounded, that is, they belong to
$L^{\infty}_{\Th}(\R^2)$ for every $\n$. The predual $\cali
M_{\n}^{\Om}(\R^2) _*$ will be denoted $A_{\n}^{\Th}(\R^2)$, and can be
interpreted as an $\n$-component of the projective Fourier algebra
$A^{\Th}(\R^2)$. In what concerns Fourier representations, the $\n$-component
$\hat\s_{\n}$ of $\l^{\Om}$ should be given by
\bea
[\hat\s_{\n}(\oml)f](q)=[(\oml\circ\k^{\n}\pd id)\cop^{\n}\hat
f_{\n}]_{\Tr_{\n}},\hs{.3} f=(\hat f_{\n})_{\Tr_{\n}}\in L^1_{\Om}\cap
L^2(\R^2).
\lb{fdecfrep}
\eea
By the same kind of manipulations made after \rep{prdeffr}, and
recalling that $A^{\Th}(\R^2)$ acts on $L^2(\R^2)$ by $\star$, the result is
$\hat\s_{\n}=id$. This should be interpreted as the injection of each
$A_{\n}^{\Th}(\R^2)$ into the von Neumann algebra $L^{\infty}_{\Th}$. The
generator of this representation can be obtained from
\bea
(g|\hat\s_{\n}(\oml_{\x\ch})\star
f)_{L^2(\R^2)}=(g\pd\ch|W^{\Th}_{\n}(f\pd\x))_{L^2(\R^2)\pd L^2(\R)}
\lb{focofrg}
\eea
and turns out to be the operator in $A_{\n}^{\Th}(\R^2)\pd\cali
M^{\Om}_{\n}(\R^2)$ given by
\bean
 W^{\Th}_{\n}(f,\x)(x,q)=e^{i\n[\Th(q';x)-\Th(q;x)]}\,f(x)\x(q+x_1),\hs 1
q\neq q'.
\eean 
The phase factors come, respectively, from the $\star$-action on $L^2(\R^2)$
($q'$ comes from the product \rep{abstpr}) and from the action of
$S^{\dagger}_{\n}(x)$ on $L^2(\R)$. From the expression \rep{wth} for the
fundamental operator $W^{\Th}$, and from the fact that the action of
$\Lo^{\dagger}(x)$ at $y$ is decomposed into the action of
$S^{\dagger}_{\n}(x)$ at $q$, we verify that $W^{\Th}_{\n}$ acts like
$W^{\Th}_{\n}\sim 1\pd S^{\dagger}_{\n}(x)$, and thus gives the genuine
decomposition of $W^{\Th}$ as the Fourier representation generator.

The coproduct and the coprojective coinvolution, when suitably decomposed from
$\K^{\Th}(\R^2)$, provide $A^{\Th}_{\n}$ with the additional structure
\bean
 \cop_{\n}f_{\n}(x,y)&=&e^{-i\n\Om(x,y)}\,f_{\n}(x+y),\\
 \k_{\n}f_{\n}(x)&=&e^{-i\n[\Th(q;x)-\Th(q;-x)]}\,f_{\n}(-x).
\eean
In the same way by which $W^{\Th}$ implements a coproduct, the generators
$W^{\Th}_{\n}$ implement the above coproducts and, as a consequence, also
satisfy the pentagonal relation.

The predual of $L^{\infty}_{\Th}$ has already been found: it is the nonabelian
algebra $L^1_{\Om}(\R^2)$. Let us examine the decomposition of its
Fourier representation $\l^{\Th}$. Since $\l^{\Th}$ is, up to a restriction on
its range of application, the left-regular representation of $L^1_{\Om}(\R^2)$,
its $\n$-component $\s_{\n}$ should be given by \rep{preweyl} with a
restriction in the range to $L^{\infty}_{\Th}(\R^2)\cap L^2(\R)$, that is,
\bea
 \s_{\n}(f)=\hat f_{\n}=\int_{\R^2} dx\, f(x)\, S_{\n}(x).
\lb{irfor}
\eea
Needless to say, these are faithful involutive representations, mapping the
twisted convolution into the projective operator product in
$\K^{\Om}_{\n}(\R^2)$, for each $\n$ in the projective dual. The generator of
this representation is easily obtained from the formula analogous to
\rep{focofrg} and is given by 
\bean
 W^{\Om}_{\n}(\x,f)(q,y)=e^{i\n[\Th(q';-y)-\Th(q;-y)]}\,\x(q-y_1)f(y),\hs 1
q'\neq q.
\eean
The same arguments which have led us to recognize $W^{\Th}_{\n}$ as the
decomposition of $W^{\Th}$ also lead to identify $W^{\Om}_{\n}$ as the
irreducible decomposition of the dual $W^{\Om}$, for they behave like
$W^{\Om}_{\n}\sim S_{\n}(y)\pd 1^o$. Furthermore, they also implement the
coproducts \rep{komncop} and consequently satisfy the pentagonal relation.

Here ends our description of the projective duality decomposition.

\section{Weyl Quantization and Duality}\lb{wqd}

We are now in condition to reexamine the Weyl-Wigner formalism in the context
of the projective Fourier duality decomposition obtained in the last section.
The expression for the projective Fourier transform \rep{preweyl} is formally
equal to the expression \rep{eldekomn} for the elements of
$\K^{\Om}_{\n}(\R^2)$, themselves given by the components $\s_{\n}$ of the
Fourier representation $\l^{\Th}$. It brings naturally to the mind Weyl's
formula, which associates a function on phase space to an irreducible
projective operator on configuration space. Before proceeding to make of this
an identification, we observe that, instead of the label $\n$, Weyl's formula
exhibits the Planck constant $\hbar$ \cite[IV,\ \S 14]{wey}. This fact leads
us to consider a rescaling in the projective dual $\Z-\{0\}$ to
$\hbar^{-1}\Z-\{0\}$, and to fix the value of the label $\n$ as $\n=1$. Doing
that means that we are selecting just one irreducible projective
representation of $\R^2$, and just one projective Kac algebra,
$\K^{\Om}_{\hbar}(\R^2)$. This shows how Quantum Mechanics is restricted to a
particular inequivalent representation, or superselection sector \cite{land}.
In this context, Weyl's formula
\bea
 \hat f_{\hbar}=\int_{\R^2}dx\,f(x)\,S_{\hbar}(x)
\lb{quasiweyl}
\eea
is a particular irreducible representation of $L^1_{\Om}$ in the operator
algebra $\K^{\Om}_{\hbar}(\R^2)$. The correspondence is completed when we write
$f$ in terms of this kind of operators. This follows from formula
\rep{wigrec}, which recover $f$ from \rep{quasiweyl} through
\bea
 f(x)=\Tr_{\hbar}[S_{\hbar}^{\dagger}(x)\hat f_{\hbar}].
\lb{recfwe}
\eea

It is also possible to rewrite Weyl's formula as a linear combination of
self-adjoint operators. This can be done by introducing operators
$\tilde S_{\hbar}(y)$ such that the projective operators $S_{\hbar}(x)$ are 
their Fourier transforms:
\bea
 S_{\hbar}(x)=\frac 1{2\p\hbar}\int_{\R^2}dy\,\ovl{\ch_x(y)}\tilde
S_{\hbar}(y),
\lb{opaad}
\eea
where $\ch_x(y)=e^{\frac i{\hbar}xy}$. Comparing $S_{\hbar}^{\dagger}(x)$ and
$S_{\hbar}(-x)$, we conclude that $\tilde S_{\hbar}^{\dagger}=\tilde
S_{\hbar}$. When we substitute \rep{opaad} in \rep{quasiweyl}, we must also
substitute the Fourier transform $\tilde f_{\hbar}$ for $f=f_{\hbar}\in
L^1_{\Om}$, 
\bean
 f(x)=\frac 1{2\p\hbar}\int_{\R^2}dz\,\ch_x(z)\,\tilde f_{\hbar}(z),
\eean
so that the two additional integrals are cancelled out by the 
characters completeness relation
\bean
\int_{\R^2} dx\,\ch_{x}(z)\ovl{\ch_{x}(z')}=(2\p\hbar)^2\d(z-z').
\eean
The Weyl formula becomes 
\bean
 \hat f_{\hbar}=\int_{\R^2}dx\,\tilde f_{\hbar}(x)\,\tilde S_{\hbar}(x),
\eean
while the Fourier transform of \rep{recfwe} gives us back the function
\bean
 \tilde f_{\hbar}(x)=\Tr_{\hbar}[\tilde S_{\hbar}(x)\hat f_{\hbar}].
\eean 
As $\tilde S_{\hbar}$ is self-adjoint, this function is real. It is the
Wigner {\em distribution function} associated to the operator $\hat f_{\hbar}$.

A particular distribution function is the Fourier transform of the
$A^{\Th}_{\hbar}(\R^2)$-function associated to the wavefunction $\x$, which is
given by \rep{prewig} with $\ch=\x$ and $\n=\hbar^{-1}$. Changing variables in
the integral, that formula can be rewritten as
\bea
\om^{\hbar}_{\x\x}(x)&=&(\x|S_{\hbar}(x)\x)\non\\
 &=&\int_{\R}dq \,e^{\frac i{\hbar}qx_2}\,\x(q+x_1/2)\ovl{\x(q-x_1/2)}.
\lb{diswig}
\eea
The Fourier transform of this function is just the Wigner distribution 
associated to the density operator $|\x\ket\bra\x|$ \cite{hil},
\bean
W_{\x}^{\hbar}(x)=[\cali F\om^{\hbar}_{\x\x}](x)=\int_{\R}dq\,e^{-\frac
i{\hbar}qx_1}\x(x_2+q/2)\ovl{\x(x_2-q/2)}.
\eean
Notice also that, if we match our notation with Dirac's through
$(\ch|S\x)=\ovl{\bra\ch|S|\x\ket}$, the function in \rep{diswig} is the same
$L^1_{\Om}$-function corresponding to the operator $\hat
f_{\hbar}=|\x\ket\bra\x|$, which is given by \rep{recfwe},
\bean
 f^{\x}(x)=\bra\x|S^{\dagger}_{\hbar}(x)|\x\ket.
\eean

The Wigner functions $\tilde f$ are also called ``quantum'' functions, since
they depend on the constant $\hbar$. Since they are the Fourier
transforms of $f\in L^1_{\Om}$, they obey the noncommutative ``twisted''
product $\circ^{\hbar}$, the Fourier image of the twisted convolution
\bean
(f\tc_{\hbar}g)(x)=\Tr_{\hbar}[S^{\dagger}_{\hbar}(x)\hat f_{\hbar}\cdot\hat
g_{\hbar}],
\eean
which is explicitly given by
\bea
[\cali F (f\tc_{\hbar}g)](z)&=&\bra\ovl{\ch_z},f\tc_{\hbar} g\ket\non\\
&=&\frac 1{2\p\hbar}\int_{\R^2\times\R^2}dx\,dy\,e^{\frac
i{\hbar}\Om(x,y)}\,\ovl{\ch_z(x+y)}f(x)g(y)\non\\ &\equiv&([\cali F
f]\circ^{\hbar}[\cali F g])(z).
\lb{tfipt}
\eea
The algebra of these essentially bounded functions with the product
$\circ^{\hbar}$ can be called the {\em Moyal algebra} and will be denoted here
by $A^{\Om}_{\hbar}(\R^2)$. The Fourier transform, which is a well known
isomorphism between the Abelian algebras $L^1$ and $L^{\infty}$ over the plane,
by \rep{tfipt} turns out to be also an isomorphism between the
noncommutative algebras $L^1_{\Om}(\R^2)$ and $A^{\Om}_{\hbar}(\R^2)$
\cite{fol,grava}. This isomorphism extends to $\cali M^{\Om}_{\hbar}$ through
$M^1_{\Om}$, for the generators $S_{\hbar}(x)$ are mapped (by \rep{recfwe})
into the ``densities'' $\d_x\in M^1_{\Om}$ and, by the Fourier transform, into
the characters $\f_x\equiv\frac 1{2\p\hbar}\ovl{\ch_x}$, the generators of the
algebra $A^{\Om}_{\hbar}$. We call $\f_x$ the generators of that algebra
because its elements are given by (the Fourier transforms)
\bean
 \tilde f_{\hbar}(y)=\int_{\R^2}dx\ \f_x(y)\,f(x)=[\cali Ff](y).
\eean
Their product, according to \rep{tfipt}, is given by
\bean
 (\f_x\circ^{\hbar}\f_y)(z)=e^{\frac i{\hbar}\Om(x,y)}\,\f_{x+y}(z),
\eean
which proves the isomorphism of $\cali M^{\Om}_{\hbar}$ and $A^{\Om}_{\hbar}$. 
Extending it further (again through the Fourier transforms) to the Kac
structure, we can add to the structure of $A^{\Om}_{\hbar}$ the coproduct and 
the coinvolution of $M^1_{\Om}$ (see \rep{estkel}):
\bean
[(\cali F\pd\cali F)\copo(f)](x,y)&=&\bra\ovl{\ch_x}\pd\ovl{\ch_y},\copo
f\ket\\ &=&\frac 1{(2\p\hbar)^2}\int_{\R^2}dz\,e^{\frac
i{\hbar}\Th(q;-z)}\,\ovl{\ch_z(x+y)}f(z)\\ &\equiv&[\cop^{\hbar}\cali
Ff](x,y);\\
\cali F\ko(f)(x)&=&\bra\ovl{\ch_x},\ko(f)\ket\\
&=&\frac 1{2\p\hbar}\int_{\R^2}dy\,e^{\frac
i{\hbar}[\Th(q;y)-\Th(q;-y)]}\,\ch_x(y)f(y)\\ &\equiv&\k^{\hbar}(\cali
Ff)(x).
\eean
Furthermore, the involution is mapped into the complex conjugation and the unit
into the constant function $\f_0=\frac 1{2\p\hbar}$, while the Haar trace
compatible with this structure is given by $\Tr^{\hbar}(\tilde
f_{\hbar})=\tilde f_{\hbar}(0)$. Summing up, the Kac structure of
$A^{\Om}_{\hbar}(\R^2)$ is given by
\bean
 \f_x\circ^{\hbar}\f_y&=&e^{\frac i{\hbar}\Om(x,y)}\,\f_{x+y};\\
 \1=\f_0&=&\frac 1{2\p\hbar};\\
 \cop^{\hbar}\f_x&=&e^{\frac i{\hbar}\Th(q;-x)}\,\f_x\pd\f_x;\\
 \k^{\hbar}\f_x&=&e^{\frac i{\hbar}[\Th(q;-x)-\Th(q;x)]}\,\ovl{\f_x};\\
 \Tr^{\hbar}(\f_x)&=&\d_x.
\eean

\section{Final Remarks}

To study how the Weyl-Wigner formalism inserts itself in the framework of
general Harmonic Analysis, we have reviewed the role of the Heisenberg and
the translation groups in the process of quantization on Euclidean phase
space. Starting from a well-established (Fourier) duality for the Heisenberg
group in terms of Kac algebras, we were able to introduce two new {\em
projective Kac algebras}, in terms of which a {\em projective duality} for
the translation group is defined. For these algebras to provide a projective
duality, the usual coinvolution axioms have being suitably adapted to the
projective framework, and this has forced us to introduce new operations. The
irreducible decomposition of the symmetric projective Kac algebra according
to the $\Om$-projective unitary dual of $\R^2$ was also performed, and it was
shown how duality survives at the irreducible level. The preduality relations
between whole and decomposed projective Kac algebras provide an explanation
for the origin of the Weyl formula as an irreducible component of the Fourier
representation of the Abelian projective Kac algebra. They also show the dual
role played by the Weyl operators and respective quantum functions, where the
latter are obtained from the first by Wigner's recovering formula and the
Fourier transform. All these facts allow us to conclude that the Weyl-Wigner
correspondence is incorporated in the projective (Fourier) duality of the
translation group as long as in the Pontryagin duality. We can go further and
ask whether it is possible to generalize this duality principle to
quantization on any other phase space. This question is partially answered
in Ref.~\cite{rasa}, where the authors have shown how far it is possible to
extend this principle to the half-plane, whose canonical group, though
requiring no central extension, has the awkward properties of being neither
Abelian nor unimodular.

In the effort towards a general quantization prescription much has yet to be
done. We have nevertheless, in the hard process of unraveling its pattern
through case-study and abstraction, got a glimpse of the basic frame and are in
position to risk a provisional proposal. Given a phase space, we should look
for its linear canonical group. Find then its two Kac algebras, the symmetric
and the Abelian. Examine the cohomology to see whether an extension is
necessary, and proceed or not to it accordingly. The resulting symmetric
algebra will be the space of quantum operators of the system.

\section*{Acknowledgments}

The authors would like to thank Prof. M. Enock for useful comments on Kac
algebras intricacies. They also would like to thank CAPES and CNPq for
financial support.

\section*{Note Added in Proof}

The authors would like to thank Prof. R. L. Hudson for pointing out a work
by him \cite{hud} where he obtained similar results, although not departing from Kac algebras. It is also worth mention that recently some related works also based on Kac algebras have been appeared \cite{enva,lara}. Although they deal with projective representations, their approach and objective are quite different from ours.

\end{document}